\documentclass[11pt,a4paper]{amsart}
\usepackage[utf8]{inputenc}
\usepackage[english]{babel}
\usepackage{amsmath}
\usepackage{amsfonts}
\usepackage{amssymb}
\usepackage{lmodern}
\usepackage[left=3.5cm,right=3.5cm,top=2.5cm,bottom=2.5cm]{geometry}

\usepackage{fancyhdr}
\usepackage{mathtools}
\usepackage{subfig,graphicx, booktabs}
\usepackage[T1]{fontenc}
\usepackage[protrusion=true, expansion=true]{microtype}
\usepackage[english]{babel}
\usepackage{xcolor}
\usepackage{tikz}
\usetikzlibrary{matrix}
\usepackage[hidelinks]{hyperref}

\usepackage[utf8]{inputenc}
\usepackage{amsmath,amssymb}
\usepackage{enumerate}
\usepackage{lmodern}

\newcommand{\rr}{\mathbb R}
\newcommand{\nn}{\mathbb N}
\newcommand{\X}{\mathcal X}
\newcommand{\Y}{\mathcal Y}
\newcommand{\e}{\mathrm{e}}

\newcommand{\dd}{\,\mathrm{d}}
\newcommand{\fcluster}{f_{\text{uni}}}
\newcommand{\vtk}{v_t^k}
\newcommand{\vtl}{v_t^{\ell}}
\usepackage{xspace}

\DeclareMathOperator*{\diag}{diag}
\DeclareMathOperator*{\law}{law}

\DeclareMathOperator*{\argmin}{argmin}

\usepackage[ruled,vlined]{algorithm2e}
\usepackage{algorithmic}

\newtheorem{theorem}{Theorem}
\newtheorem{lemma}{Lemma}
\newtheorem{corollary}{Corollary}
\theoremstyle{plain}

\newtheorem{remark}{Remark}

\title[Consensus-Based Optimization for Multi-Objective Problems]{Consensus-Based Optimization for Multi-Objective Problems: A Multi-Swarm Approach}
\author[K.~Klamroth, M.~Stiglmayr, C.~Totzeck]{Kathrin Klamroth \and Michael Stiglmayr \and Claudia Totzeck}
\address{University of Wuppertal, Gau{\ss}str. 20, 42119 Wuppertal}
\email{\{klamroth,stiglmayr,totzeck\}@uni-wuppertal.de}

\begin{document}

\begin{abstract}We propose a multi-swarm approach to approximate the Pareto front of general multi-objective optimization problems that is based on the Consensus-based Optimization method (CBO). The algorithm is motivated step by step beginning with a simple extension of CBO based on fixed scalarization weights. To overcome the issue of choosing the weights we propose an adaptive weight strategy in the second modelling step. The modelling process is concluded with the incorporation of a penalty strategy that avoids clusters along the Pareto front and a diffusion term that prevents collapsing swarms. Altogether the proposed $K$-swarm CBO algorithm is tailored for a diverse approximation of the Pareto front and, simultaneously,  the efficient set of general non-convex multi-objective problems. The feasibility of the approach is justified by analytic results, including convergence proofs, and  a performance comparison to the well-known non-dominated sorting genetic algorithm (NSGA2) and the recently proposed one-swarm approach for multi-objective problems involving Consensus-based Optimization.
\end{abstract}

\maketitle

\section{Introduction}

Multiple conflicting objective functions occur in a variety of applications ranging from engineering design to economic and financial decisions. Economical goals are often in conflict with ecological criteria, we have to trade-off between expected return and risk, and we aim at affordable yet high quality products. We refer to the textbooks \cite{Ehrg05} and \cite{Miettinen1999} for a general introduction to the topic of multi-objective optimization. 
In this paper, we focus on \emph{continuous multi-objective optimization problems} (MOP)
\begin{equation}\tag{MOP}\label{eq:MOP}
	\min\limits_{x \in \X} f(x)=(f_1(x),\ldots,f_p(x)),
\end{equation}
with a non-empty 
feasible set $\X\subset\rr^d$, $d\in\nn$, and with $p\geq 2$ continuous objective functions $f_i\colon \X\rightarrow\rr$, each of which has a unique global minimum on $\X$. 
We consider unconstrained problems, i.e., $\X=\rr^d$, as well as box-constrained problems where $\X=[\ell,u]$ with lower and upper bounds $\ell,u\in\rr^d$ with $\ell_j\leq u_j$, $j=1,\dots,d$.

We refer to $\rr^d$ as the \emph{decision space} and to $\rr^p$ as the \emph{objective space} of \eqref{eq:MOP}. The set $\Y\coloneqq f(\X)\subset\rr^p$ is called the feasible \emph{outcome set} of \eqref{eq:MOP}. Two feasible solutions $x^1,x^2\in\X$ are compared based on their respective outcome vectors $z^1=f(x^1)$ and $z^2=f(x^2)$: We say that $z^1$ \emph{dominates} $z^2$ (and $x^1$ dominates $x^2$)
denoted by $z^1\leqslant z^2$ if and only if
$$ 
z^1_i\leq z^2_i \text{~for all~} i=1,\dots,p \text{~and~} z^1\neq z^2.
$$

A feasible solution $x\in\X$ is called \emph{Pareto optimal} or \emph{efficient} if there is no other solution $\bar{x}\in\X$ such that $f(\bar{x})\leqslant f(x)$. The corresponding image in the objective space is called \emph{non-dominated} in this case. The set of all Pareto optimal solutions is denoted by $\X_P$ and referred to as the \emph{efficient set}, and the set of all non-dominated outcome vectors, i.e., $\Y_P=f(\X_P)$, is called the \emph{non-dominated set} or the \emph{Pareto front} of \eqref{eq:MOP}. Note that since we assume that each objective function has a unique global minimum on $\X$, 
the Pareto front of \eqref{eq:MOP} is non-empty since it contains these individual minima. 

An important approach to generate or to approximate Pareto optimal solutions are scalarizations that transform the multi-objective problem \eqref{eq:MOP} into 
a series of associated single-objective problems. Maybe the most prominent scalarization approach is the \emph{weighted-sum scalarization} \cite{gass:thec:1955}: Given non-negative weights $\lambda_i\geq0$, $i=1,\ldots,p$ (that correspond to the relative importance of the respective criteria), the \emph{weighted sum-objective} is given by
\begin{equation}\label{eq:weightedsum}
	f_{\lambda}(x) \coloneqq \sum_{i=1}^p \lambda_i\, f_i(x).
\end{equation}
Let $\Lambda\coloneqq\{\lambda\in\rr^p \colon \sum_{i=1}^p\lambda_i=1 \text{~and~} \lambda_i>0,\, i=1,\dots,p\}$ and let $\Lambda_0\coloneqq\{\lambda\in\rr^p \colon \sum_{i=1}^p\lambda_i=1 \text{~and~} \lambda_i\geq0,\, i=1,\dots,p\}$. Further, we define $\rr^p_+\coloneqq\{z\in\rr^p\colon z_i\geq 0,\; i=1,\ldots,p\}$. The following theorem is a well-known result from the field of multi-objective optimization.
\begin{theorem}[see, e.g., \cite{Ehrg05}] 
	\label{thm:efficient_set}
	If $\lambda\in\Lambda$, then an optimal solution $\bar{x}(\lambda)\in\X$ of \eqref{eq:weightedsum} is efficient for \eqref{eq:MOP}. Moreover, if $\Y+\rr^p_+$ is convex and $\bar{x}\in\X$ is efficient, then there exists a weighting vector $\bar{\lambda}\in\Lambda_0$ such that $\bar{x}=\bar{x}(\bar{\lambda})$ is optimal for \eqref{eq:weightedsum}.
\end{theorem}
For later reference we phrase the following observations as remark.
\begin{remark}[see, e.g., \cite{Ehrg05}]\label{rem:Thm1}
	When $\lambda\in\Lambda_0$, then an optimal solution $\bar{x}(\lambda)\in\X$ of \eqref{eq:weightedsum} is at least \emph{weakly efficient} for \eqref{eq:MOP}, i.e., there is no $\hat{x}\in\X$ such that $f_i(\hat{x})<f_i(\bar{x}(\lambda))$ for all $i=1,\ldots,p$. Note also that when $\Y+\rr^p_+$ is non-convex, then it can not be guaranteed in general that all efficient solutions can be obtained as an optimal solution of a weighted-sum scalarization. 
\end{remark}

The goal of this paper is to develop a provably convergent, yet efficient algorithm for high-quality representations of Pareto fronts of multi-objective optimization problems. We choose the \emph{Consensus-based Optimization method (CBO)} as basis for the algorithm as it is a  particle method for single-objective global optimization problems which is easy to implement and allows for analytical studies. Its dynamics is governed by an interacting particle system that allows to propagate information of the individuals through a weighted mean. The weighted mean has two advantages: (a) there is no need to label individuals as current best, and (b) there is no need for pairwise interactions as the dynamic of one particle depends only on the weighted mean of the whole swarm. Advantage (a) allows for the derivation of a corresponding mean-field equation \cite{CBOmeanfield} which can be employed for analytical studies such as the large time behaviour and the convergence to the global minimizer. Indeed, in \cite{Jin,CBO2} it is shown that the invariant solution of the mean-field dynamics is a Dirac-delta located arbitrarily close to the global minimizer of the objective function. 

For completeness we recall the CBO dynamics for $N$ particles in the single-objective case, i.e., for $p=1$. Hence, let $f \colon \mathbb R^d \rightarrow \mathbb R_+$ denote the objective function. For $\lambda, \sigma >0$ the dynamics of the $j$-th particle $X^j:[0,T]\rightarrow\rr^d$, $j=1,\dots,N$, at time $t\in[0,T]$ is given by
\begin{subequations}\label{eq:CBO}
	\begin{gather}
		\dd X_t^{j} = -(X_t^{j} - v_t)\, \dd t + \sigma \,\diag(X_t^{j} - v_t) \dd B_t^{j}, \\ \law(X_0^{j}) = \rho_0, \qquad j=1,\dots, N,
	\end{gather}   
	where the weighted mean $v_t$ is given by
	\begin{equation}
		v_t = \frac{\sum_{j=1}^{N} X_t^{j} \, \e^{-\alpha f(X_t^{j})}}{\sum_{j=1}^{N}  \e^{-\alpha f(X_t^{j})}},
	\end{equation}
\end{subequations}
$\rho_0$ is a probability distribution on the state space and $B^j:[0,T]\rightarrow\rr^d$, $j=1,\dots,N$ are independent Brownian motions. Here and in the following we use the anisotropic noise term $\sigma \diag(X_t^{j} - v_t) \dd B_t^{j}$ as proposed in \cite{Jin} as it is shown to be more robust in settings with high-dimensional state space.

The idea of the multi-swarm CBO algorithm we propose in Section~\ref{sec:fixedWeights} and further develop in Section~\ref{sec:adaptiveWeights} is based on Theorem~\ref{thm:efficient_set}. Indeed, each swarm is associated with a weight vector that yields a scalarization of the cost function. Following the CBO dynamic the swarms globally minimize their respective scalarized problems giving us an approximation of the Pareto front. For a diverse approximation we introduce interactions between the swarms in Section~\ref{sec:adaptiveWeights}. The multi-swarm CBO algorithm with adaptive swarms is analyzed in Section~\ref{sec:proof_diverse_approximation}, where we show under appropriate assumptions that each swarm clusters at a point along the Pareto front and the approximation obtained by the points of all swarms is diverse in the sense that the clustering points admit a distance greater than a minimal distance.

As mentioned above, weighted sum scalarization can have issues in case of non-convex problems. In particular, efficient solutions do not necessarily correspond to global minima of the weighted sum scalarization. By introducing a penalization strategy (in Section~\ref{sec:general}) that avoids clustering along the Pareto front we allow swarms to also converge towards efficient unsupported solutions, i.e., solutions which are not located on the convex hull of \(\Y\).
In the modelling process we focus on the weighted means of the swarms and tailor the dynamics such that the means form a high quality, i.e., representative and well-distributed approximation of the Pareto front. In Consensus-based Optimization the swarms are forced to collapse in the large time limit \cite{CBO2,CBO1}. This is disadvantageous in our setting, as we would loose local information about the Pareto front. Hence, we prevent collapsing swarms by implementing the diffusion term from Consensus-based sampling \cite{CarrilloHoffmannStuartVaes}. 

Each modelling step is illustrated by simulation results. To not interrupt the flow, we present the details of the simulations in Section~\ref{sec:parameters}, where all parameters and implementation details are reported. We conduct a qualitative study comparing the proposed multi-swarm CBO with the recently proposed single-swarm CBO for multi-objective problems \cite{BorghiHertyPareschi} and the well-known NSGA2 algorithm implemented in \cite{pymoo} in Section~\ref{sec:comparison}, before we conclude the article and give an outlook to future work.

\section{CBO for convex multi-objective problems}
\label{sec:simple}
We begin with the presentation of a multi-swarm approach for multi-objective optimization based on CBO (MSCBO)
with fixed scalarization weights and provide analytical results. Fixed weights come with the burden of choosing appropriate weights a priori, which may result in highly unbalanced Pareto front approximations  \cite{DasDennis}. To overcome this problem, we propose an adaptive weight strategy in the second part of this section.

\subsection{MSCBO 
	with fixed scalarization weights}\label{sec:fixedWeights}

We begin with a simple CBO method for multi-objective problems that is gradually extended in the following: Let $K\in \nn$ be the number of swarms involved and let $(\lambda^k)_{k=1,\dots,K} \subset \Lambda$ be their respective fixed scalarization weights. For each $k$ we introduce the swarm size $N_k \in \nn$ and the positions of the individuals are denoted by $X^k = (X^{k,1},\dots,X^{k,N_k}) \colon [0,T] \rightarrow \rr^{N_k\cdot d}$. The evolution of the swarms is given by the dynamics
\begin{subequations}\label{eq:independent_dynamics}
	\begin{gather}
		\dd X_t^{k,j} = -(X_t^{k,j} - {\vtk}) \dd t + \sigma\, \diag(X_t^{k,j} - {\vtk}) \dd B_t^{k,j}, \\ \law(X_0^{k,j}) = \rho_0, \qquad j=1,\dots, N_k, \; k=1,\dots,K,
	\end{gather}   
	where the weighted mean ${\vtk}$ is given by 
	\begin{equation}\label{eq:weightedaverage}
		{\vtk} = \frac{\sum_{j=1}^{N_k} X_t^{k,j} \e^{-\alpha f_{\lambda^k}(X_t^{k,j})}}{\sum_{j=1}^{N_k}  \e^{-\alpha f_{\lambda^k}(X_t^{k,j})}},
	\end{equation}
\end{subequations}
$\sigma>0$ is a diffusive strength parameter, $\rho_0 \in \mathcal P_2(\rr^d)$ is a probability measure with finite second moment, $B^{k,j}$ are independent Brownian motions and $\alpha>0$ allows to scale the difference of the local and global minima in the objective function (Laplace principle). For notational convenience we define $N=\sum_{k=1}^K N_k$ and the vector $X_t^N = (X_t^1, \dots, X_t^K) \in \rr^{N\cdot d}.$ A pseudo code for the dynamics is given in Algorithm~\ref{alg:fixed}.

\begin{remark}[Initialization of the weights]\label{rem:weights}
	In order to obtain admissible initial conditions for the weight vectors $\lambda^k, k=1,\dots,K,$ we draw uniformly distributed random samples from the standard simplex $\mathcal S^p=\Lambda_0$. 
	Thus, each weight vector has $p-1$ degrees of freedom and we can write $\lambda^k = (\lambda^{k,1}, \dots, \lambda^{k,p-1}, 1- \sum_{i=1}^{p-1} \lambda^{k,i}).$  For the numerical tests we use Algorithm~2 from \cite{uniformSimplex} to sample the initial weights $\lambda^k$ for $k=1,\dots,K$.
\end{remark}

\begin{remark}
	In contrast to the dynamics proposed in \cite{borghi22repulsion,BorghiHertyPareschi}, where the interacting particles in the swarm converge towards different efficient solutions, here each swarm (namely its weighted mean) represents one solution.
	That means each swarm is supposed to find one point along the Pareto front and the swarms act independently. In Section~\ref{sec:adaptiveWeights} we implement interactions between the swarms through adaptive scalarization weights.
\end{remark}

\subsubsection{Analysis of the independent swarms}
Note that in the dynamics proposed above the $K$ swarms are independent and each is globally minimizing its weighted sum \eqref{eq:weightedsum} given by $\lambda^k.$ Hence, we may employ the result in \cite{CBO_MTW} to prove the well-posedness of the dynamics. 
\begin{theorem}[Well-posedness \cite{CBO_MTW}]
	Let $K\in \nn, (\lambda^k)_{k=1,\dots,K} \subset \Lambda_0, N_k \in \nn$ fixed and $f_{\lambda^k}$ locally Lipschitz continuous for every $\lambda^k$, then system~\eqref{eq:independent_dynamics} admits a unique strong solution $\{ X_t^N \colon t\ge 0 \}$ for each initial condition $X_0^{N}$ with $\mathbb E[|X_0^{N}|^2] < \infty.$
\end{theorem}
\begin{proof}
	By the independence of the swarms, we first apply Theorem 1 in \cite{CBO_MTW} to every subsystem describing the evolution of each swarm and collect the corresponding solutions in $X^N$ for the full system. 
\end{proof}

As the convergence analysis on the particle level is highly nontrivial, we follow the steps in \cite{Jin,CBO2} and derive statistical representations of the swarms in form of mean-field equations. Under the standard assumption of \textit{propagation of chaos}, i.e.~that the particles decouple as $N_k \rightarrow \infty$, we denote the probability of finding a member of the $k$-th swarm in position $x$ at time $t$ by $\rho_t^k(x).$ Following the lines of \cite{Jin,CBO2} we compute the evolution equation of these probabilities with the help of Itô-calculus.
Towards this end, let for simplicity $N_k = \bar N$ for all $k=1,\dots,K$. For $\bar N \rightarrow \infty$ the PDE system corresponding to  
\eqref{eq:independent_dynamics} reads
\begin{align*}
	\partial_t \rho_t^k &= 
	\frac{\sigma^2}{2} \sum_{i=1}^d \partial_{ii}\big((x - v[\rho_t^k] ) \,\rho_t^k\big) 
	+ \nabla \cdot \bigl( (x - v[\rho_t^k])\, \rho_t^k \bigr), \\
	\lim\limits_{t\rightarrow 0} \rho_t^k &= \rho_0, \quad 
	k=1,\dots,K, 
\end{align*}
with the weighted mean given by
\[
v[\rho_t^k] = \frac{1}{\int \e^{-\alpha f_{\lambda^k}(x)} \dd \rho_t^k(x)}\int x\, \e^{-\alpha f_{\lambda^k}(x)} \dd \rho_t^k(x)
\]
and $\sigma > 0,$ $\rho_0 \in \mathcal P_2(\rr^d)$ as above.

Again, using the independence of the $K$ swarms we can directly apply the result in \cite{CBO2} to prove the convergence of each swarm towards a limit point arbitrarily close to the efficient set. To state the result properly we define
\begin{align*}
	E^k(t) &\coloneqq \int_{\rr^d} x \dd \rho_t^k(x),\\
	V^k(t) &\coloneqq \int_{\rr^d} |x - E^k(t) | \dd \rho_t^k(x), \\
	M^k(t) &\coloneqq \int_{\rr^d} \e^{-\alpha  f_{\lambda^k}(x)} \dd \rho_t^k(x).
\end{align*}
The convergence theory for CBO \cite{CBO2} is based on the (strong) assumption that the minimizer of the objective function is unique. 
This can be guaranteed under a similarly strong uniqueness assumption for
the minimizer of the scalarized objective functions, i.e.,
\begin{equation} \label{ass:1} \tag{A1} 
	\bar{x}^k\coloneqq\argmin_{x\in\X}f_{\lambda^k}(x) \text{~is unique, and} \quad 
	F^k \coloneqq f_{\lambda^k}(\bar{x}^k)>0 \quad \text{for all $\lambda^k$} 
\end{equation}
In addition, we inherit the following technical assumption from the convergence theory of CBO \cite{CBO2}
\begin{equation} \label{ass:2} \tag{A2}
	c_f^k \coloneqq \max\bigl\{ \| \max\limits_{i} |\partial_{ii} f_{\lambda^k}| \|_\infty, \; \| r(\nabla^2 f_{\lambda^k})\|_\infty \bigr\} < \infty, 
\end{equation}
where $\nabla^2 f_{\lambda^k}$ represents the Hessian of $f_{\lambda^k},$ $r(\nabla^2 f_{\lambda^k})$ is the spectral radius and $\partial_{ii} f_{\lambda^k}$ is the $i$-th element of the diagonal of $\nabla^2 f_{\lambda^k}.$ 

\begin{remark}
	Since we assumed that all objective functions $f_i$, $i=1,\dots,p$, are bounded from below, all weighted sums $f_{\lambda^k}$ are bounded from below for $\lambda^k\in\Lambda_0$. Hence, $F^k>0$ in Assumption~\eqref{ass:1} can always be guaranteed after applying an appropriate linear transformation to the objective functions values. The uniqueness of the global minimizers of weights sum scalarizations $f_{\lambda^k}$, however, does in general not follow from the uniqueness of the global minimizers of the individual objective functions $f_i$, $i=1,\dots,p$, even if $p=2$. As an example, suppose that $\X=[0,1]\subset\rr$, $f_1(x)=x$ and $f_2(x)=1-x$. Then, $\Y=\Y_P$ is the line segment connecting the two points $(0,1)$ and $(1,0)$ in $\rr^2$, and while the global minimizers of $f_1$ and of $f_2$ are unique, the complete feasible set $\X$ is optimal for  $f_{\lambda}$ with $\lambda=(\frac{1}{2},\frac{1}{2})$.
\end{remark}

\begin{theorem}[Convergence towards the efficient set \cite{Jin}] \label{thm:approximation}
	Let \eqref{ass:1}  and \eqref{ass:2} hold. Further, let $\alpha$ and the initial condition $\rho_0$ be chosen such that for each $k \in \{1,\dots,K \}$ it holds
	\begin{equation*}
		\mu^k \coloneqq 2 - \sigma^2 - 2\, \sigma^2 \,\frac{\e^{-\alpha F^k}}{M^k(0)} >0, \qquad 
		\nu^k \coloneqq \frac{2\, V^k(0)}{\mu^k\cdot (M^k(0))^2} \,\alpha\, \e^{-2\alpha F^k} c_f^k (2 + \sigma^2) \leq \frac{3}{4}.
	\end{equation*}
	Then for each $k \in \{1,\dots,K\}$ it holds that $V^k(t) \rightarrow 0$ exponentially fast and there exist $\tilde x^k$ such that $v[\rho_t^k] \rightarrow \tilde x^k$, $E^k(t) \rightarrow \tilde x^k$ exponentially fast. Moreover, for $\alpha \rightarrow \infty$ the limit points $\tilde x^k$ are arbitrarily close to the efficient set. In particular, for every $\epsilon > 0$ there exists $\alpha \gg 1$ such that $\tilde x^k \in B_\epsilon(\bar{x}^k).$ 
\end{theorem}
\begin{proof}
	As the $k$ swarms are independent the result is a direct consequence of Theorem~3.1 in \cite{Jin}.
\end{proof}

\begin{corollary}
	Let $f_i$ be strictly convex functions for all $i\in \{1,\dots,p\}$ and let Assumption~\eqref{ass:2} be satisfied. Then for each point of the efficient set there exists $\lambda^k \in \Lambda_0$ such that the swarm corresponding to $\lambda^k$ and following the dynamics \eqref{eq:independent_dynamics} concentrates arbitrarily close to a point in the efficient set.
\end{corollary}
\begin{proof}
	The result follows from Theorem~\ref{thm:efficient_set} in combination with Theorem~\ref{thm:approximation}.
\end{proof}

The previous results are based on the fact that the weights are fixed. It therefore remains open to choose the weights appropriately. This is not a trivial task, as it is for example well-known that an equidistant choice of weights does not necessarily lead to an equidistant resolution of the Pareto front, even in the case of convex problems \cite{DasDennis}.
In order to circumvent the problem of choosing appropriate weights, we propose an adaptive procedure including dynamic weights in the next section.

\subsection{MSCBO with dynamic weight adaption}\label{sec:adaptiveWeights}
Already in the biobjective case the choice of scalarization weights is nontrivial. In fact, assuming that $f_1$ and $f_2$ are continuously differentiable, the optimal choice of the weights depends on the ratio $\frac{\partial f_2}{\partial f_1}$ \cite{DasDennis}, which is in general unknown. 

This observation motivates to incorporate adaptive weights in the multi-objective CBO method. The adaption dynamics are written as an ODE. In order to circumvent the restriction $\lambda^k_j \in [0,1],$ we use the bijective transformation 
\[
\text{for } \lambda^k \in \Lambda \text{ set } \mu^k= \ln(\lambda^k)  \text{ elementwise}
\]
and formulate the dynamic weight adaption in terms of $\mu^k.$
The scalarized cost function then reads
\[
f_{\mu^k}(x) = \frac{\sum_{i=1}^p \exp(\mu_i^k) f_i(x)}{\sum_{i=1}^p \exp(\mu_i^k) }.
\]

\begin{remark}
	Note that $\lambda^k \in \Lambda$ instead of $\Lambda_0$ implies that only efficient solutions with bounded trade-off can be determined. Since efficient solutions with unbounded trade-off usually correspond to extreme cases that are not very interesting from a practical point of view, this does not impose a strong restriction.
\end{remark}

The dynamic of adaptive weights is based on a pairwise interaction of swarms \(k\) and \(\ell\) if their respective weight vectors $\mu^{k}$ and $\mu^{\ell}$ are close to each other, or the distance of  their weighted means, \(\vtk\) and \(\vtl\), in the objective space is small, i.e., if the swarms map to similar points \( f(\vtk)\) and \(f(\vtl)\) in the objective space. Let
\[ 
d_{k,\ell} \coloneqq|\mu^{k} - \mu^{\ell}|, \qquad d_{k,\ell}^f \coloneqq | f(\vtk) - f(\vtl) | < c_\text{rep} \leq \infty.
\]
Then the repulsion is modelled by means of interaction potentials, where in the case $c_\text{rep} < \infty$ these potentials are assumed to have compact support. 

For the numerical tests we use Morse potentials \cite{Dorsogna} which are given by
\begin{equation}\label{eq:potential}
	P(d_{k,\ell}) = R\, \e^{-d_{k,\ell}/r} - A\, \e^{-d_{k,\ell}/a}, \quad R,A \geq 0,\; r,a > 0.
\end{equation}
The parameters $R, A$ denote the strength of the repulsion and attraction, respectively, while $r, a$ allow to adjust the range of the attractive and repulsive force, respectively. 

\begin{remark}
	We choose to have strong repulsive forces on a very short range and attractive forces on a medium range. Hence on larger ranges interactions can be neglected and we consider only interactions with direct neighbors in our main theorem analysing the approximation of the Pareto front (Theorem~\ref{thm:diversity}).
\end{remark}

For interacting particle systems, see for example \cite{Schafe1}, the force on $\mu^{k}$ that results from the interaction with $\mu^\ell \neq\mu^k$  
is given by the corresponding gradient
\begin{equation}\label{eq:potentialforce}
	\nabla_{\mu^{k}} P(d_{k,\ell}) = \left(\frac{A}{a} \, \e^{-d_{k,\ell}/a} - \frac{R}{r} \, \e^{-d_{k,\ell}/r} \right) \frac{\mu^{k} - \mu^{\ell}}{|\mu^{k} - \mu^{\ell}|}.
\end{equation}
The main goal we are pursuing with the adaptive weights is to obtain a diverse approximation of the Pareto front. We therefore add another interaction strategy that depends on the distance of the weighted means in the objective space. Following the spirit of \eqref{eq:potentialforce} we replace the distances $d_{k,\ell}$ in the prefactor by the distance of the objective function values $d_{k,\ell}^f,$ leading to the second force given by

\begin{equation}\label{eq:forceF}
	\left(\frac{A_f}{a_f} \, \e^{-d_{k,\ell}^f/a_f} - \frac{R_f}{r_f} \, \e^{-d_{k,\ell}^f/r_f} \right) \frac{\mu^{k} - \mu^{\ell}}{|\mu^{k} - \mu^{\ell}|} , \quad R_f,A_f \geq 0,\; r_f,a_f > 0.
\end{equation} 
Altogether, this leads to the following forces for the adaptive weight adjustment 
\begin{align*}
	&\mathcal K(X_t^k,X_t^\ell,\mu^k, \mu^\ell) =  \\
	&=\left(\frac{A}{a} \, \e^{-d_{k,\ell}/a} - \frac{R}{r} \, \e^{-d_{k,\ell}/r} 
	+ \frac{A_f}{a_f} \, \e^{-d_{k,\ell}^f/a_f} - \frac{R_f}{r_f} \, \e^{-d_{k,\ell}^f/r_f} \right) 
	\frac{\mu^k - \mu^\ell}{|\mu^k - \mu^\ell|}.
\end{align*} 
More generally, one can use any force with similar properties to define an adaptive weight adjustment given by
\begin{align}\label{eq:adative_weights}
	\frac{\dd}{\dd t} \,\mu^k &= - \frac{1}{\tau}\sum_{\ell=1, \ell\ne k}^K \mathcal K(X_t^k,X_t^\ell,\mu^k,\mu^\ell), \\
	\mu^k(0) &= \ln(\lambda_0^k), \quad k=1,\dots,K. \notag
\end{align}
Here $\tau>0$ is a parameter that allows us to control the time scale of the interaction dynamics. In fact, it will become important in Section~\ref{sec:proof_diverse_approximation}, where we assume that the adaption of the scalarization weights happens much faster than the interaction dynamics and concentration process of the swarms.

We recover the original weights 
\[
\lambda^k = \frac{\exp(\mu^k)}{\sum_k \exp(\mu^k)}.
\]

Combining \eqref{eq:independent_dynamics} and  \eqref{eq:adative_weights} yields a CBO dynamic for multi-objective problems with self-adaptive scalarization weights. We want to emphasize that this combination introduces interdependencies between swarms. Moreover, the dynamics of the weights lead to dynamic scalarized cost functions. A pseudo code of this method can be found in Algorithm~\ref{alg:adaptive}.

\subsubsection{Numerical comparison of static and dynamic weights}\label{sec:results_static_dynamic}
\begin{figure}[ht!]
	\centering
	\subfloat[\label{fig:static_dynamics_weightsA}]{\includegraphics[scale=0.35,clip,trim=45 0 35 30] {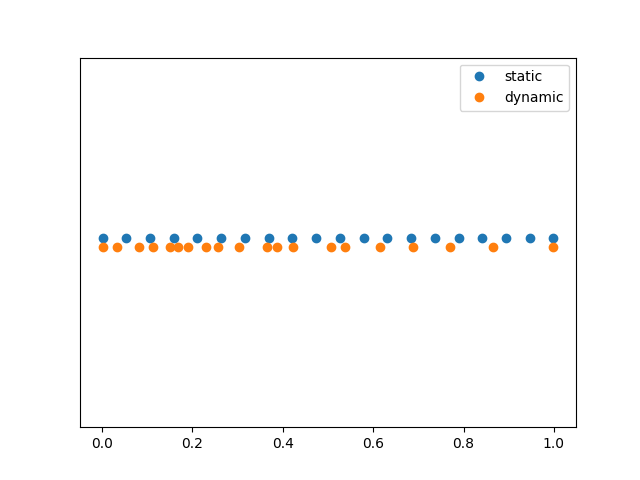}}
	\hspace{\fill}
	\subfloat[\label{fig:static_dynamics_weightsB}]{\includegraphics[scale=0.35,clip,trim=10 0 35 30] {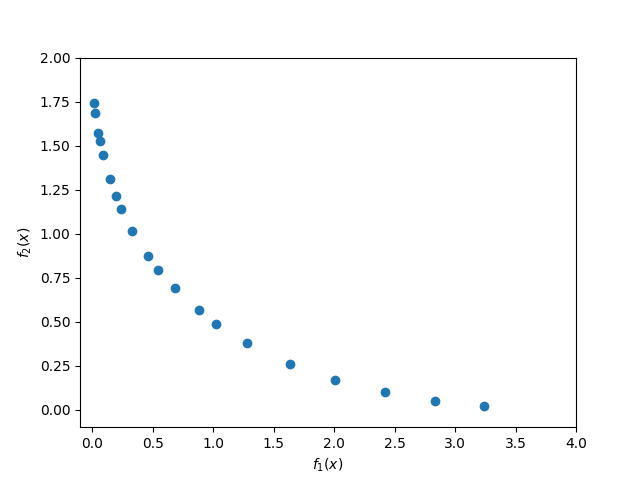}}
	\hspace{\fill}
	\subfloat[\label{fig:static_dynamics_weightsC}]{\includegraphics[scale=0.35,clip,trim=10 0 35 30] {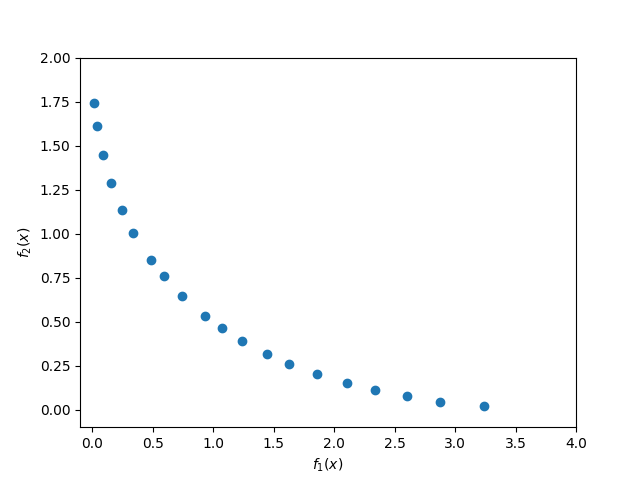}}
	\caption{Illustration of adaptive weights. The orange points in the plot on the left show the adaptive weights at time $t=2.5$. The corresponding front given by $f({\vtk})_{k=1,\dots,10}$ is shown in the plot on the right. In the middle the approximation obtained with equally distributed weights is shown.}
	\label{fig:static_dynamics_weights}
\end{figure}

To illustrate the effect of the dynamic weight adaption, we solve the test problem \textsf{Schaffer1}, see Section~\ref{sec:implementation} for more details.
The weight vectors are initialized equidistantly in $(0,1)$ 
like shown with the blue points in Figure~\ref{fig:static_dynamics_weightsA}. The corresponding front found by solving \eqref{eq:independent_dynamics} is shown in Figure~\ref{fig:static_dynamics_weightsB}. The distribution of the adaptive $\lambda^k$ is given by the orange points in  Figure~\ref{fig:static_dynamics_weightsA} and the corresponding front in Figure~\ref{fig:static_dynamics_weightsC}. We note that equally distributed weights stress the left part of the front. The adaptive weighting counteracts this effect by shifting the weights to the right, which results in a better resolution of the front. All other details regarding the numerical results are given in Section~\ref{sec:implementation}.

The results indicate that the proposed method works fine for convex problems. However, we note that we may run into problems in non-convex settings. Indeed, as discussed in Remark~\ref{rem:Thm1}, we can only find convex parts of the front by globally minimizing weighted sum scalarizations. For general and possible non-convex problems we propose another extension in the next section. Before we enter the part of the general problems we analyse the method with adaptive weights in terms of well-posedness, mean-field approximation and convergence.

\subsubsection{Well-posedness of the SDE system with adaptive weights}
To analyse the SDE system with adaptive weights we introduce the notation $Y_t = [X_t, \mu_t]$ where $X_t^N = (X_t^k)_{k=1,\dots,K}$
and $\mu_t^N = (\mu_t^k)_{k=1,\dots,K}.$ Moreover, we write the SDE system in vectorized form using the notation
\begin{subequations}\label{eq:adaptiveSystem}
	\begin{gather}
		\dd Y_t^N = -G(Y_t^N) \dd t + \sigma D(Y_t^N) \dd B_t^N, \\
		G(Y_t^N) = \bigl[G_X (Y_t^N), G_{\mu}(Y_t^N)\bigr], \qquad D(Y_t^N) = \diag(G_X(Y_t^N), \textbf{0}), \\
		Y_0^N = [X_0^N, \mu_0^N].
	\end{gather}
\end{subequations}
where 
\begin{align*}
	\big(G_X(Y_t^N)\big)_k &= \bigl( (X_t^{k,i} - {\vtk}) \bigr)_{i=1,\dots,N_k}, \text{ and} \\
	\big(G_\mu(Y_t^N)\big)_k &= \frac{1}{\tau} \Bigl(\sum_{\ell=1, \ell \ne k}^K \mathcal K(X_t^k, X_t^\ell,\mu^k, \mu^\ell)\Bigr). 
\end{align*}
To obtain a well-posedness result for the system with adaptive weights, we make the following assumptions:

The interaction force $\mathcal K$ is locally Lipschitz and grows at most linearly in $x,y,\mu$ or 
$\nu$. In particular, it holds
\begin{equation}\label{ass:locLip}\tag{A3}
	\bigl|\mathcal K(x,y,\mu,\nu) - \mathcal K(\bar x, \bar y, \bar \mu, \bar \nu)\bigr| \leq C_L \bigl(|x-\bar x| + |y-\bar y| + |\mu-\bar \mu| + |\nu-\bar \nu|\bigr)
\end{equation}
and 
\begin{equation}\label{ass:linGrowth}\tag{A4}
	\bigl|\mathcal K(x,y,\mu,\nu)\bigr| \leq C_b \bigl(|x| + |y| + |\mu| + |\nu|\bigr)
\end{equation}
for some constants $C_L, C_b >0$ and $|(x,y,\mu,\nu)-(\bar x,\bar y, \bar \mu, \bar \nu)|<M.$

\begin{remark}\label{rem:smoothPotential}
	As usual in particle interaction dynamics the forces resulting from physical potentials such as the Morse potential proposed above, do not satisfy the regularity assumption of ODE theory. In order to comply with these assumptions we have to consider a smoothed version, as for example,
	$$\hat P(d) = R^{-d^2/r} - Ae^{-d^2/a}, \quad R,A\ge 0,\quad r,a >0. $$
\end{remark}

Similar to the single objective case, the proof of the well-posedness of the SDE system with adaptive weights is based on the following technical lemma.
\begin{lemma}\label{lem:technical}
	Let $K\in\nn, N_k \in \nn$ for all $k=1,\dots,K$, $N = \sum_{k=1}^K N_k,$ further $\alpha,M > 0$ and \eqref{ass:locLip}, \eqref{ass:linGrowth} hold. Then there exist constants $C_1,C_2 >0$ depending only on $M$ such that for any $Y_t = [X_t,\mu_t], \hat Y_t =[\hat X_t, \hat \mu_t] \in \rr^{d N} \times \rr^{pK}$ 
	with $|Y_t|, |\hat Y_t| \leq M$ it holds
	\begin{align*}
		\bigl|\big(G_X(Y_t) - G_X(\hat Y_t) \big)_{k,j}\bigl| 
		+ \bigl|\big(G_\mu(Y_t) - G_\mu(\hat Y_t) \big)_{k}\bigr|   &\leq C_1 |Y_t^k - \hat Y_t^k|, \\ 
		|G_X(Y_t)| + |G_\mu(Y_t)| &\leq C_2 |Y_t^k|.
	\end{align*}
\end{lemma}
\begin{proof}
	Let $N = \sum_{k=1}^K N_k,$ $Y_t = [X_t,\mu_t], \hat Y_t =[\hat X_t, \hat \mu_t] \in \rr^{dN} \times \rr^{pK}$. 
	We begin with
	\begin{align*}
		&\big(G_X(Y_t) - G_X(\hat Y_t) \big)_{k,j} =\\
		&= \frac{\sum_{m=1}^{N_k} (X_t^{k,j} - X_t^{k,m})\, \e^{-\alpha f_{\mu^k}(X_t^{k,m})}}{\sum_{m=1}^{N_k} \e^{-\alpha f_{\mu^k}(X_t^{k,m})}} - \frac{\sum_{m=1}^{N_k} (\hat X_t^{k,j} - \hat X_t^{k,m}) \,\e^{-\alpha f_{\mu^k}(\hat X_t^{k,m})}}{\sum_{m=1}^{N_k} \e^{-\alpha f_{\mu^k}(\hat X_t^{k,m})}} \\
		&\quad+ \frac{\sum_{m=1}^{N_k} (\hat X_t^{k,j} - \hat X_t^{k,m})\, \e^{-\alpha f_{\mu^k}(\hat X_t^{k,m})}}{\sum_{m=1}^{N_k} \e^{-\alpha f_{\mu^k}(\hat X_t^{k,m})}} - \frac{\sum_{m=1}^{N_k} (\hat X_t^{k,j} - \hat X_t^{k,m})\, \e^{-\alpha f_{\hat \mu^k}(\hat X_t^{k,m})}}{\sum_{m=1}^{N_k} \e^{-\alpha f_{\hat \mu^k}(\hat X_t^{k,m})}} \\
		&=: I_1 + I_2 + I_3 + I_4 + I_5
	\end{align*}
	where
	\begin{align*}
		I_1 &= \frac{\sum_{m=1}^{N_k} \big(X_t^{k,j} - \hat X_t^{k,j} - (X_t^{k,m} - \hat X_t^{k,m}) \big) \,\e^{-\alpha f_{\mu^k}(X_t^{k,m})}}{\sum_{m=1}^{N_k} \e^{-\alpha f_{\mu^k}(X_t^{k,m})}}, \\
		I_2 &= \frac{\sum_{m=1}^{N_k} (\hat X_t^{k,j} - \hat X_t^{k,m})  (\e^{-\alpha f_{\mu^k}(X_t^{k,m})} - \e^{-\alpha f_{\mu^k}(\hat X_t^{k,m})})}{\sum_{m=1}^{N_k} \e^{-\alpha f_{\mu^k}(X_t^{k,m})}}, \\
		I_3 &=\sum_{m=1}^{N_k} (\hat X_t^{k,j} - \hat X_t^{k,m})  \e^{-\alpha f_{\mu^k}(\hat X_t^{k,m})} \frac{\sum_{m=1}^{N_k} (\e^{-\alpha f_{\mu^k}(\hat X_t^{k,m})} - \e^{-\alpha f_{\mu^k}( X_t^{k,m})})}{ \sum_{m=1}^{N_k} \e^{-\alpha f_{\mu^k}( X_t^{k,m})} \sum_{m=1}^{N_k} \e^{-\alpha f_{\mu^k}(\hat X_t^{k,m})} }, \\
		I_4 &=  \frac{\sum_{m=1}^{N_k} (\hat X_t^{k,j} - \hat X_t^{k,m})  (\e^{-\alpha f_{\mu^k}(\hat X_t^{k,m})} - \e^{-\alpha f_{\hat \mu^k}(\hat X_t^{k,m})})}{\sum_{m=1}^{N_k} \e^{-\alpha f_{\mu^k}(\hat X_t^{k,m})}}, \\
		I_5 &= \sum_{m=1}^{N_k} (\hat X_t^{k,j} - \hat X_t^{k,m})\, \e^{-\alpha f_{\hat \mu^k}(\hat X_t^{k,m})} \frac{\sum_{m=1}^{N_k} (\e^{-\alpha f_{\mu^k}(\hat X_t^{k,m})} - \e^{-\alpha f_{\hat\mu^k}( \hat X_t^{k,m})})}{ \sum_{m=1}^{N_k} \e^{-\alpha f_{\mu^k}( \hat X_t^{k,m})} \sum_{m=1}^{N_k} \e^{-\alpha f_{\hat \mu^k}(\hat X_t^{k,m})} }.
	\end{align*}
	Similar to \cite{CBO2} the terms can be bounded as follows:
	\begin{align*}
		|I_1| &\leq |X_t^{k,j} - \hat X_t^{k,j}| + |X_t^k - \hat X_t^k|, \\
		|I_2| &\leq \sqrt{\frac{2}{N_k}} \, c_X^k \, |X_t^{k} - \hat X_t^k| \, \sqrt{N_k |\hat X_t^{k,j}|^2 + |\hat X_t^k|^2}, \\
		|I_3| &\leq \frac{\sqrt{2}\, c_X^k}{N_k} \, |X_t^k - \hat X_t^k| \,\sqrt{|\hat X_t^{k,j}|^2 + |\hat X_t^k|^2},
	\end{align*} 
	where $c_X^k \coloneqq \alpha \| \nabla f_{\mu^k} \|_{L^\infty(B_X)} \, \exp(\alpha \| f_{\mu^k} - \underline{ f_{\mu^k}} \|_{L^\infty(B_X)})$ for $B_X = \{ x \in \rr^d \colon |x| \leq M \}$ and $\underline{ f_{\mu^k}} \coloneqq \min\limits_{x \in B_X} f_{\mu^k}(x).$
	To bound $I_4$ and $I_5$ we first note that
	\begin{align*}
		f_{\mu^k} (x) - f_{\hat \mu^k}(x) = \frac{\sum_{i=1}^p (\e^{\mu_i^k} - \e^{\hat \mu_i^k}) \, f_i(x)}{\sum_{i=1}^p \e^{\mu_i^k}} 
		+ \sum_{i=1}^p \e^{\hat\mu_i^k} f_i(x) \frac{\sum_{i=1}^p (\e^{\hat\mu_i^k} - \e^{\mu_i^k})}{\sum_{i=1}^p \e^{\hat\mu_i^k} \, \sum_{i=1}^p \e^{\mu_i^k}}
	\end{align*} 
	which leads to the estimate
	\begin{equation*}
		|f_{\mu^k} (x) - f_{\hat \mu^k}(x)|  \leq \Bigl(1+\frac{1}{p}\Bigr) \, \| f \|_{L^\infty(B_X)} \,\e^{2M} \,|\mu^k - \hat\mu^k|   
	\end{equation*}
	for $|\mu^k|, |\hat \mu^k| \leq M.$ Using this we obtain
	\begin{align*}
		|I_4| &\leq \sqrt{\frac{2}{N_k}} \, d_M^k \, \sqrt{N_k |\hat X_t^{k,j}|^2 + |\hat X_t^k|^2} \, \| f \|_{L^\infty(B_X)}\, \Bigl(1+\frac{1}{p}\Bigr) \, \e^{2M} \,|\mu^k - \hat\mu^k|, \\
		|I_5| &\leq \frac{\sqrt{2}}{N_k} \, d_M^k \, \sqrt{|\hat X_t^{k,j}|^2 + |\hat X_t^k|^2} \, \| f \|_{L^\infty(B_X)}\,  \Bigl(1+\frac{1}{p}\Bigr) \, \e^{2M} \, |\mu^k - \hat\mu^k|, 
	\end{align*} 
	where $d_M^k \coloneqq \alpha \,\exp(2\alpha \| f \|_{L^\infty(B_x)}).$
	As these estimates are independent of $j=1,\dots,N_k$ and $k=1,\dots,K$, this proves that $G_X$ is locally Lipschitz. Note that $G_\mu$ is locally Lipschitz by Assumption~\eqref{ass:locLip}. 
	
	Moreover, it is easy to see that $|(G_X(Y_t))_{k,j}| \leq |X_t^{k,j}| + |X_t^k| \leq 2 |Y_t^k|$ holds which leads to the linear growth of $G_X.$ Further, the linear growth of $G_\mu$ is ensured by Assumption~\eqref{ass:locLip}. 
\end{proof}

Based on this Lemma we establish the following well-posedness and existence result.
\begin{theorem}\label{thm:wellposednessAdaptive}
	Let the assumptions of Lemma~\ref{lem:technical} be satisfied. For each $K\in\nn$ and $N_k \in \nn$ for $k=1,\dots,K,$ system \eqref{eq:adaptiveSystem} admits a unique strong solution $\{Y_t \colon t\ge 0\}$ for any initial condition $Y_0$ satisfying $\mathbb{E}[|Y_0|^2] < \infty.$
\end{theorem}
\begin{proof}
	The proof is based on \cite[Chapter 5.3, Theorem~3.1]{Durrett}. Note that due to Lemma~\ref{lem:technical} we only have to show that there exists $b >0$ such that it holds 
	\[
	-2 \,Y_t \cdot G(Y_t) + \sigma^2\, \text{trace}\bigl(D(Y_t)D(Y_t)^\top \bigr) \leq b \, |Y_t|^2.
	\]
	For the first term we obtain, using \eqref{ass:locLip} and \eqref{ass:linGrowth},
	\[
	-2 \, Y_t \cdot G(Y_t) \leq 2\, C_1 \,|Y_t|^2
	\]
	for some constant $C_1>0.$
	Moreover, by the diagonal structure of $D(Y_t)$ we get
	\[
	\text{trace}\bigl(D(Y_t)\, D(Y_t)^\top \bigr) \leq 2\, C_2 \, |Y_t|^2
	\]
	for some constant $C_2>0.$ Altogether, this proves the desired estimate and  Theorem~3.1 in \cite[Chapter 5.3]{Durrett} yields the result.
\end{proof}

\begin{remark}
	Note that the regularized potential proposed in Remark~\ref{rem:smoothPotential} satisfies the regularity assumptions of Theorem~\ref{thm:wellposednessAdaptive}.
\end{remark}

\subsubsection{Mean-field limit with adaptive weights}
Formally passing to the mean-field limit yields a coupled PDE/ODE system given by
\begin{align*}
	\partial_t \rho_t^k &= \frac{\sigma^2}{2} \sum_{i=1}^d \partial_{ii}\big((x - v[\rho_t^k] ) \rho_t\big) + \nabla \cdot \big((x - v[\rho_t^k]\big) \rho_t^k), \\
	\lim\limits_{t\rightarrow 0} \rho_t^k &= \rho_0, \quad k=1,\dots,K,  \\
	\frac{\dd}{\dd t} \mu^k &= -\frac{1}{\tau}\sum_{\ell=1, \ell\ne k}^K \mathcal K(\rho_t^k,\rho_t^\ell, \mu^k, \mu^\ell), \qquad
	\mu^k(0) = \ln(\lambda_0^k)
\end{align*}
with the weighted mean given by
\[
v[\rho_t^k] = \frac{1}{\int \e^{-\alpha f_{\lambda^k}(x)} d\rho_t^k(x)}\int x\, \e^{-\alpha f_{\lambda^k}(x)} \dd\rho_t^k(x), \qquad \lambda^k = \exp(\mu^k)
\]
and $\tau,\sigma > 0,$ $\rho_0 \in \mathcal P_2(\rr^d)$ as above. Note that the  coupling through the weights makes the proof of convergence towards the efficient set considerably more difficult as there is a trade-off between the concentration of the swarms and the balancing of the weight vectors. This is addressed in the next section.

\subsection{Diverse approximation of the Pareto front}\label{sec:proof_diverse_approximation}
The aim of this section is to prove that the proposed scheme leads to a diverse approximation of the Pareto front in the case of strictly convex biobjective problems, i.e.~$f \colon \rr^{d} \rightarrow \rr^2.$ In fact, we prove that each swarm concentrates at one point along the Pareto front in the long-time limit, $t\rightarrow \infty.$ The proof builds on results of the single-objective CBO scheme in \cite{CBO2,fornasier2021globally}. Certainly, we require additional assumptions that we motivate and formulate in the following.

For the sake of a simple presentation of the idea, we consider the smoothed version of the Morse interaction potential discussed in  Remark~\ref{rem:smoothPotential}. Moreover, we assume that the initial condition of the scalarization weights and the parameters of the interaction potentials are chosen such that the scalarization weights do not leave the domain $(\epsilon_\lambda, 1-\epsilon_\lambda)\times(\epsilon_\lambda, 1-\epsilon_\lambda)$ for fixed $0<\epsilon_\lambda\ll1.$  This can be ensured by choosing $A_f,R_f,a_f$ and $r_f$ such that we have short-range repulsion and mid-range attraction and negligible interaction in the long-range. Note that this assumption allows us to work directly with normalized scalarization weights $\lambda^k = (\lambda_1^k,\lambda_2^k),\; \lambda^k_1 + \lambda_2^k=1$ for all $k=1,\dots,K.$  In fact we do not require the reformulation in terms of $\mu$ that is used for the numerics. 

To be more precise, we consider the dynamics
$$\frac{\dd}{\dd t}\lambda^k = \frac{-1}{\tau}\sum_{\ell=1}^K \left(\frac{A_f}{a_f} \, \e^{-(d_{k,\ell}^f)^2/a_f} - \frac{R_f}{r_f} \, \e^{-(d_{k,\ell}^f)^2/r_f} \right) (\lambda^k - \lambda^\ell) \quad \text{for}\; k=1,\dots,K.$$

Note that 
\[\lambda^k - \lambda^\ell = \begin{pmatrix} \lambda^k_1 - \lambda^\ell_1 \\ \lambda^\ell_1 -\lambda^k_1 \end{pmatrix}=\begin{pmatrix} \lambda^\ell_2 - \lambda^k_2 \\ \lambda^k_2 -\lambda^\ell_2 \end{pmatrix} \quad \text{for all $k,\ell$}.
\] 
As all other quantities on the right-hand side are scalar, the $\lambda^k$ stay on the constraint $\lambda_1^k + \lambda_2^k=1$ for all times $t>0$ once they are normalized.

As a consequence, the dynamics of the vectors $\lambda^k$ are uniquely defined by the dynamics of one of their components. Without loss of generality we chose the first component here and obtain
\begin{equation*} \frac{\dd}{\dd t} \lambda^k_1 = -\sum_{\ell=1}^K u(d_{k,\ell}^f) (\lambda^k_1 -  \lambda^{\ell}_1), \qquad \frac{\dd}{\dd t} \lambda^{k}_{2} = -\sum_{\ell=1}^K u(d_{k,\ell}^f) ( \lambda^\ell_1 - \lambda^k_1) = 1-\frac{\dd}{\dd t} \lambda^k_1,
\end{equation*}
where we use the notation
$$ u(d_{k,\ell}^f) = \left(\frac{2 A_f}{a_f} \, \e^{-(d_{k,\ell}^f)^2/a_f} - \frac{2R_f}{r_f} \, \e^{-(d_{k,\ell}^f)^2/r_f} \right). $$
Note that the unique representation by the first component allows us to sort the scalarization weights of the swarms. Without loss of generality, we assume in the following that $\lambda_1^1 < \lambda_1^2 < \dots < \lambda_1^K.$ Moreover, due to negligible long-range interaction we only consider interactions with direct neighbors.

The strategy of the proof is as follows: first, we exploit the relationship of the non-dominated points $f(\bar x(\lambda))$ corresponding to a given $\lambda\in\Lambda$ (i.e., $\bar x(\lambda)$ is the unique optimum of \eqref{eq:weightedsum} with weight $\lambda$) and the scalarization weight $\lambda$ itself. In fact, we assume that 
\begin{eqnarray*}
	\begin{pmatrix} \lambda^k_1-\lambda^\ell_1\\ \lambda^\ell_1-\lambda^k_1\end{pmatrix}
	&=& S(\lambda^k, \lambda^\ell) 
	f(\bar{x}(\lambda^k))-f(\bar{x}(\lambda^\ell))\quad\text{for all $k,\ell$}
\end{eqnarray*}
for some problem dependent matrix $S(\lambda^k, \lambda^\ell).$
Second, using the chain rule, we notice that
$$\frac{\dd}{\dd t} f(\bar x(\lambda^k)) = \frac{\dd}{\dd x} f(\bar x(\lambda^k)) \frac{\dd}{\dd \lambda_1^k} \bar x(\lambda^k) \frac{\dd}{\dd t} \lambda_1^k. $$
Now, we employ a Lyapunov argument for all swarms at once, which will show that  the non-dominated points $f(\bar x(\lambda^k))$ spread along the Pareto front with $\|\cdot\|_2$-distance $d_\text{min}$  in the long time limit, where $d_\text{min}$ is the unique root of $u(d).$  

Let us state the assumptions and the theorem: 
\begin{enumerate}[({A}1)]\setcounter{enumi}{4}
	\item The Pareto front of the biobjective problem $f\colon \rr^{d} \rightarrow \rr^2$ is strictly convex and connected.
	\item  The parameters of the interaction potential satisfy $A,R=0$ and $A_f,R_f,a_f$ and $r_f$ are such that the interaction potential $\mathcal U$ models short-range repulsion and mid-range attractions; long-range interactions are negligible. In particular, we assume that there is a unique root $d_\text{min}$ of $u(d)$ such that all $K$ non-dominated points $f(\bar x(\lambda^k)), k=1,\dots K$ can spread along the front with Euclidean distance greater than $d_\text{min}.$ Moreover, we assume that each swarm is only interacting with its direct neighbors (in the sense of the sorted weight vectors).
	\item The initial values of the scalarization weights $\lambda^k$ satisfy $\lambda^k_1 + \lambda^k_2 = 1$ and $(\lambda^k_1, \lambda^k_2) \in [\epsilon_\lambda, 1-\epsilon_\lambda]^2$ for each $k=1,\dots,K$ and fixed $0 <\epsilon_\lambda \ll 1.$ Moreover, $\lambda_1^1 < \lambda_1^2 < \dots < \lambda_1^K$ and the corresponding non-dominated points $f(\bar x(\lambda^k))$ are distributed such that $\epsilon_\lambda \le \lambda_1 ^1$ and $\lambda_1^K \le 1-\epsilon_\lambda.$ 
	
	\item For each pair $\lambda^k, \lambda^\ell$ with $k,\ell\in\{1,\dots,K\}$, $k\neq\ell$, there exists a negative definite matrix $S(\lambda^k, \lambda^\ell)\in\rr^{2\times 2}$ such that 
	\[
	\begin{pmatrix}
		\lambda^k_1 - \lambda^\ell_1 \\  \lambda^\ell_1 - \lambda^k_1
	\end{pmatrix} = S(\lambda^k, \lambda^\ell) (f(\bar x(\lambda^k))- f(\bar x (\lambda^\ell))) \text{ and } S(\lambda^k, \lambda^\ell) = S(\lambda^\ell, \lambda^k).
	\]
	\item For each $\lambda \in [\epsilon_\lambda, 1- \epsilon_\lambda]$ there exists a negative definite matrix $T(\lambda)$ such that
	\[
	\frac{\dd}{\dd t} f(\bar x(\lambda)) = T(\lambda) \frac{\dd}{\dd t}\lambda = T(\lambda) \frac{\dd}{\dd t} \begin{pmatrix} \lambda_1 \\ 1-\lambda_1 \end{pmatrix}. 
	\]
	\item The product $\mathcal T(\lambda) \mathcal S(\lambda)$ with 
	{\small \begin{equation*} \mathcal T(\lambda) = 
			\begin{tikzpicture}[baseline=(current bounding box.center)]
				\matrix (m) [matrix of math nodes,nodes in empty cells,right delimiter={)},left delimiter={(},inner sep=-2pt,nodes={inner sep=1ex}]{
					T(\lambda^1) + T(\lambda^2) & -T(\lambda^2)& & & \\
					-T(\lambda^2) & T(\lambda^2) + T(\lambda^3) & - T(\lambda^3) & & \\
					& - T(\lambda^3) &  &  &\\
					& & & & -T(\lambda^{K-1})  \\
					& & & -T(\lambda^{K-1}) & T(\lambda^{K-1}) + T(\lambda^K)\\
				} ;
				\draw[loosely dotted,thick] (m-2-2)-- (m-5-5);
				\draw[loosely dotted,thick] (m-2-3)-- (m-4-5);
				\draw[loosely dotted,thick] (m-3-2)-- (m-5-4);
			\end{tikzpicture}
	\end{equation*}}
	and
	\begin{equation}
		\mathcal S(\lambda) = \begin{pmatrix} S(\lambda^1,\lambda^2) & & \\ & \ddots \\ & & S(\lambda^{K-1},\lambda^K) \end{pmatrix}
	\end{equation} 
	is symmetric. Here $T(\lambda^k)$ are as defined in (A9) and $S(\lambda^k,\lambda^\ell)$ as defined in (A8) for each $k,\ell=1,\dots,K$, $k\neq\ell$.
	\item  The dynamics of the scalarization weights and the dynamics of the swarms have different time scales. In fact, the adaption of the scalarization weights is much faster than the dynamics of the swarms. 
\end{enumerate}
We first state the result for the limiting case $\alpha = \infty$ as this ensures that $f(v[\rho_t^k]) = f(\bar x(\lambda))$ by the Laplace principle \cite{fornasier2021globally}.
\begin{theorem}\label{thm:diversity}
	Let (A1)-(A11) hold and $\alpha = \infty.$ For $K\in\nn$ the dynamics 
	\begin{align*}
		\partial_t \rho_t^k &= \frac{\sigma^2}{2} \sum_{i=1}^d \partial_{ii}\big((x - v[\rho_t^k] ) \rho_t\big) + \nabla \cdot \big((x - v[\rho_t^k]\big) \rho_t^k), \quad \lim\limits_{t\rightarrow 0} \rho_t^k = \rho_0,\\
		\frac{\dd}{\dd t}\lambda^k &= -\frac{1}{\tau}\sum_{\ell=1}^K \left(\frac{A_f}{a_f} \, \e^{-(d_{k,\ell}^f)^2/a_f} - \frac{R_f}{r_f} \, \e^{-(d_{k,\ell}^f)^2/r_f} \right) (\lambda^k - \lambda^\ell), \\
		v[\rho_t^k] &= \frac{1}{\int \e^{-\alpha f_{\lambda^k}(x)} d\rho_t^k(x)}\int x\, \e^{-\alpha f_{\lambda^k}(x)} \dd\rho_t^k(x), \quad k=1,\dots,K, 
	\end{align*}
	yields a diverse approximation of the Pareto front of $f$ in the long time limit. In particular, each swarm concentrates at $\bar x(\lambda^k)$  with pairwise distance $$\| f(\bar x(\lambda^k)) - f(\bar x(\lambda^\ell))\|_2 \ge d_\text{min} \text{ for all } k\ne \ell \text{ as } t\rightarrow \infty.$$
	In other words, the stationary solution of the dynamics consists of $K$ Dirac-measures located at the  points $\bar x(\lambda^k)$  for $k=1,\dots,K$ which have a pairwise  Euclidean distance greater or equal $d_{\min}.$
\end{theorem}
\begin{proof}
	First note that due to $\alpha = \infty,$ we have have $v_f[\rho_t^k] = \bar x(\lambda^k)$ for $t>0$ \cite{fornasier2021globally}.  This allows us to focus only on the dynamics of the scalarization weights in the following. Indeed, we follow the steps motivated above and begin with some observations on the pairwise interactions. Using (A6), (A8) and (A9) we obtain 
	\begin{align*}
		\frac{\dd}{\dd t} f(\bar x(\lambda^k)) &= \frac{\dd}{\dd x} f(\bar x(\lambda^k)) \frac{\dd}{\dd \lambda^k} \bar x(\lambda^k) \frac{\dd}{\dd t} \lambda^k = T(\lambda^k) \frac{\dd}{\dd t} \lambda^k \\
		&=  -T(\lambda^k) \sum_{\ell=k-1}^{k+1} u(d_{k,\ell}^f) (\lambda^k - \lambda^\ell) \\
		&= -T(\lambda^k) \sum_{\ell=k-1}^{k+1}  S(\lambda^k,\lambda^\ell) u(d_{k,\ell}^f)(f(\bar x(\lambda^k)) - f(\bar x(\lambda^\ell)))
	\end{align*}
	for $k=1,\dots,K$. For simplicity we set the summand to zero whenever $\ell < 1$ or $\ell > K.$
	
	In our strictly convex setting, the sorting of the weights $\lambda^k$ induces a sorting of the points $f(\bar x(\lambda^1)), \dots, f(\bar x(\lambda^K)$ in decending order long the $f_1$-axis.
	To reconstruct the vectors $f(\bar x(\lambda^k)), k=1,\dots,K$ at $t>0$, it is sufficient to know their initial values and the relative positions from $k$ to $k+1$. This will be exploited in the following.
	
	For $k=1,\dots K-1$ we define $q_{k,k+1} = f(\bar x(\lambda^k)) - f(\bar x(\lambda^{k+1})).$ Clearly, it holds
	\begin{align*} \frac{\dd}{\dd t} q_{k,k+1} &= \frac{\dd}{\dd t} f(\bar x(\lambda^k)) - \frac{\dd}{\dd t} f(\bar x(\lambda^{k+1})) \\
		&= -T(\lambda^k) \sum_{\ell=k-1}^{k+1}  S(\lambda^k,\lambda^\ell) u(d_{k,\ell}^f)(f(\bar x(\lambda^k)) - f(\bar x(\lambda^\ell))) \\
		&\qquad +T(\lambda^{k+1}) \sum_{\ell=k}^{k+2}  S(\lambda^{k+1},\lambda^\ell) u(d_{k+1,\ell}^f)(f(\bar x(\lambda^{k+1})) - f(\bar x(\lambda^\ell))).
	\end{align*}
	Replacing the differences yields
	\begin{align*} 
		&\frac{\dd}{\dd t} q_{1,2} = -\big(T(\lambda^1) + T(\lambda^2)\big) S(\lambda^1,\lambda^2) u(d_{1,2}^f) q_{1,2} \\
		&\qquad \qquad + T(\lambda^2) S(\lambda^2,\lambda^3)u(d_{2,3}^f) q_{2,3} \\
		&\frac{\dd}{\dd t} q_{K-1,K} = T(\lambda^{K-1})S(\lambda^{K-2},\lambda^{K-1}) u(d_{K-2,K-1}^f) q_{K-2,K-1} \\
		&\quad \qquad \qquad -\big(T(\lambda^{K-1}) +T(\lambda^{K})\big) S(\lambda^{K-1},\lambda^{K}) u(d_{K-1,K}^f) q_{K-1,K}  
	\end{align*}
	and for all $k=2,\dots, K-2$ it holds
	\begin{align*} 
		&\frac{\dd}{\dd t} q_{k,k+1} =  T(\lambda^{k})S(\lambda^{k-1},\lambda^{k}) u(d_{k-1,k}^f) q_{k-1,k} \\
		&\quad \qquad \qquad -\big(T(\lambda^{k}) + T(\lambda^{k+1}) \big)S(\lambda^{k},\lambda^{k+1}) u(d_{k,k+1}^f) q_{k,k+1}  \\
		&\quad \qquad \qquad + T(\lambda^{k+1}) S(\lambda^{k+1},\lambda^{k+2})u(d_{k+1,k+2}^f) q_{k+1,k+2}
	\end{align*}

	Now, let us consider the vector $q = (q_{1,2}, \dots, q_{K-1,K})$ and the potential $$\mathcal V (q) = \sum_{k=1}^{K-1} \left( R_f e^{-|q_{k,k+1}|^2/r_f} - A_fe^{-|q_{k,k+1}|^2/a_f}\right)  $$ 
	with gradient
	$$ \nabla \mathcal V(q) = \left( u(d_{1,2}^f) q_{1,2}, \dots, u(d_{K-1,K^f}) q_{K-1,K} \right)^\top.$$
	The vectorized dynamics reads
	\begin{equation*}
		\frac{\dd}{\dd t} q = -\mathcal T(\lambda) \mathcal S(\lambda) \nabla \mathcal V(q),
	\end{equation*}
	with $\mathcal T(\lambda)$ and $\mathcal S(\lambda)$ as in (A10).
	
	It is easy to check using the binomial formulas and (A9) that $\mathcal T(\lambda)$ is negative definite. By (A10) $\mathcal T(\lambda)\mathcal S(\lambda)$ is symmetric and therefore positive definite. By (A7) we can bound the smallest eigenvalue of the product, $\sigma_0(\lambda)>0$ by some constant depending on $\epsilon_\lambda$, i.e., $\sigma_0(\lambda) > c_{\epsilon_\lambda} > 0$ for all $\lambda.$
	
	The positive definiteness of $\mathcal T(\lambda)\mathcal S(\lambda)$ allows us to define a norm $\| \cdot \|_{\mathcal T\mathcal S}$ induced by the scalar product $\langle \cdot, T(\lambda)\mathcal S(\lambda) \cdot \rangle$. 
	The time evolution of $\mathcal V(q)$ can now be computed as
	\begin{align}
		\frac{\dd}{\dd t} \mathcal V(q) =  -\| \nabla \mathcal V(q) \|_{\mathcal T\mathcal S}^2 \le -c_{\epsilon_\lambda} \| \nabla V(q) \|^2 < 0.
	\end{align}
	Hence $\mathcal V$ is a Lyapunov functional for $q$ and hence in the long-time limit the dynamics stabilizes such that $\nabla \mathcal V$ = 0. In particular, by (A6) that means the Euclidean distances of the non-dominated points is greater than $d_{\min}.$ 
	
	Now, $\tau$ in (A11) allows us to balance the times of the stabilization of the scalarization weights and the time to collapse the swarms. For $\tau \rightarrow 0$ the scalarization weights converge very fast to a diverse approximation of the Pareto front. The scalarization weights are then stationary and the swarms concentrate at their weighted average as discussed in Theorem~\ref{thm:approximation}. This concludes the proof.
\end{proof}

We want to emphasize that in this setting each of the swarms admits a stationary solution which is the efficient point along the Pareto front where the swarm concentrates. In contrast, the particles of the dynamic proposed in \cite{BorghiHertyPareschi} may move along the front for all times as Theorem~4.1 in \cite{BorghiHertyPareschi} shows that the density of the particles move towards the Pareto front, but no stationarity is shown. 
Before presenting an example satisfying (A5)-(A10) we comment on the case $\alpha < \infty.$
\begin{remark}
	We expect to obtain a similar result for $\alpha < \infty$ using the quantitative Laplace principle \cite{fornasier2021globally}. Clearly, the points $f(v[\rho_t^k]), k=1,\dots,K$ will only be in an $\epsilon$-neighborhood of the Pareto front, where $0<\epsilon\ll1$ depends on $\alpha.$ The details of the proof are beyond the scope of this article. 
\end{remark}

\subsubsection{Example for (A5)-(A10)}
A priori it is not clear that any problem satisfies the assumptions (A5)-(A10). We therefore present an example with two swarms represented by the weighting vectors $\lambda^1,\lambda^2\in\Lambda$ in the following. Note that biobjective convex quadratic optimization problems have been extensively studied in the  literature, therefore parts of the following analysis, including optimality conditions and the subsequent sensitivity analysis, can also be found, e.g., in \cite{art:BOLT21,hill:nonl:2001,tour:onbi:2019}. 

Let us consider the biobjective stricly convex quadratic problem 
$$ f :\rr\rightarrow\rr^2, \qquad  f(x)=\begin{pmatrix} x^2 +1\\ \frac12 (x-1)^2 +1 \end{pmatrix}. $$ 
Let $\lambda\in[0,1]$ and consider the scalarized objective
$$ f_\lambda(x) = \lambda (x^2+1) + (1-\lambda) \big(\frac12(x-1)^2 + 1\big).$$
Note that this is equivalent to a weighted sum scalarization \eqref{eq:weightedsum} with $(\lambda_1,\lambda_2)=(\lambda,1-\lambda)$. 
The efficient solutions parameterized by $\lambda\in[0,1]$ can be computed using the first order optimality condition for $f_\lambda$ w.r.t.\ the variable $x\in\rr$: 
\begin{align*}
	\frac{\dd}{\dd x} f_\lambda(x) &= (1+\lambda)x - (1-\lambda) \stackrel{!}{=} 0 \\
	\implies& \bar{x}(\lambda) = \frac{1-\lambda}{1 + \lambda} \quad\text{and}\quad f(\bar{x}(\lambda))=\begin{pmatrix} \left(\frac{1-\lambda}{1+\lambda}\right)^2 +1 \\ 2\left(\frac{\lambda}{1+\lambda}\right)^2 +1 \end{pmatrix}. 
\end{align*}
The function $P \colon [0,1] \rightarrow \rr^2$ given by $P(\lambda) = f(\bar x(\lambda))$ is continuous, therefore (A5) holds. 
Moreover, we obtain for the difference of the non-dominated points of the two swarms as
\begin{eqnarray*}
	f(\bar{x}(\lambda_1^{1}))-f(\bar{x}(\lambda_1^2))&=&
	\begin{pmatrix} \left(\frac{1-\lambda_1^1}{1+\lambda_1^1}\right)^2-\left(\frac{1-\lambda_1^2}{1+\lambda_1^2}\right)^2 \\ 2\left(\frac{\lambda_1^1}{1+\lambda_1^1}\right)^2-2\left(\frac{\lambda_1^2}{1+\lambda_1^2}\right)^2 \end{pmatrix}\\
	&=&\frac{-2}{(1+\lambda_1^1)(1+\lambda_1^2)}\begin{pmatrix} \frac{1-\lambda_1^1}{1+\lambda_1^1} + \frac{1-\lambda_1^2}{1+\lambda_1^2} & 0 \\ 0 & \frac{\lambda_1^1}{1+\lambda_1^1} + \frac{\lambda_1^2}{1+\lambda_1^2}\end{pmatrix}
	(\lambda^1 - \lambda^2).
\end{eqnarray*}
As we assume that all entries of $\lambda^1, \lambda^2$ are bounded away from zero and one, the matrix is invertible. This allows us to define $S(\lambda^1, \lambda^2)$ through its inverse
$$S(\lambda^1, \lambda^2)^{-1} = \frac{-2}{(1+\lambda_1^1)(1+\lambda_1^2)}\begin{pmatrix} \frac{1-\lambda_1^1}{1+\lambda_1^1} + \frac{1-\lambda_1^2}{1+\lambda_1^2} & 0 \\ 0 & \frac{\lambda_1^1}{1+\lambda_1^1} + \frac{\lambda_1^2}{1+\lambda_1^2}\end{pmatrix}.$$
Concerning (A9) we compute
\begin{align*} \frac{\dd}{\dd t} f(\bar x(\lambda^1)) &= \frac{-2}{(1+\lambda_1^1)^3} \begin{pmatrix}
		0 & 2 \\ -2 & 0 
	\end{pmatrix} \, \begin{pmatrix}\lambda_1^1 \\ 1- \lambda_1^1 \end{pmatrix} \frac{\dd}{\dd t} \lambda_1^1  \\
	&= \frac{-2\, u(d_{1,2}^f)}{(1+\lambda_1^1)^3} \begin{pmatrix} 2(1-\lambda_1^1) & 0 \\ 0 & 2 \lambda_1^1 \end{pmatrix}  (\lambda^1 - \lambda^2) 
\end{align*}
Hence, we find $$ T(\lambda^1) =  \frac{-2\, }{(1+\lambda_1^1)^3} \begin{pmatrix} 2(1-\lambda_1^1) & 0 \\ 0 & 2 \lambda_1^1 \end{pmatrix}$$
is negative definite and analogous for $k=2.$ As $T(\lambda^1), T(\lambda^2)$ and $S(\lambda^1,\lambda^2)$ are all diagonal, it is easy to check that (A10) holds.

\medskip
Up to here, we mainly focused on convex objective functions since only convex parts of the efficient set can be obtained by global minimization of weighted sum scalarizations, see \cite{DasDennis} for more details. 
In Section~\ref{sec:general} we introduce a penalization of clusters in the objective space in order to obtain well-dispersed representations of the Pareto front and, as an intended  side-effect, reach into non-convex parts of the Pareto front.

\section{MSCBO for general multi-objective problems}
\label{sec:general}
As discussed in the previous section, we cannot expect the present dynamics to work well in non-convex settings. We therefore propose a penalization strategy in order to obtain dynamics that lead to reasonable approximations for general multi-objective problems.

\subsection{Penalization of clusters in the objective space}
To motivate the penalization strategy we analyze the approximation obtained with the adaptive weight algorithm applied to a test problem with \textsf{dent}, see Section~\ref{sec:testproblems} for more details.  

\begin{figure}[ht!]
	\centering
	\subfloat[\label{fig:A_interactionFA}]{\includegraphics[scale=0.38]{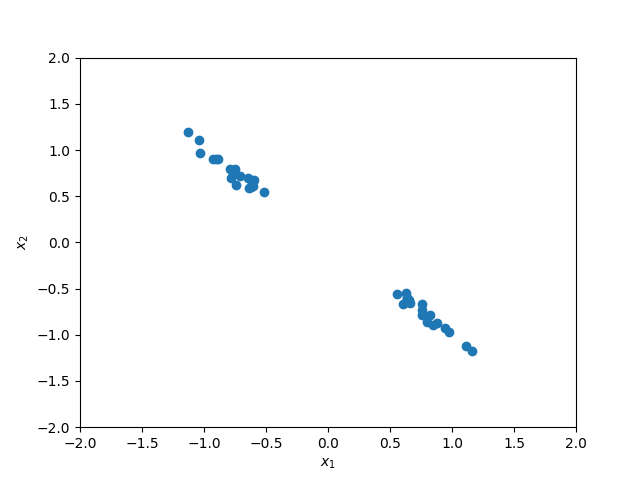}}
	\subfloat[\label{fig:A_interactionFB}]{\includegraphics[scale=0.38]{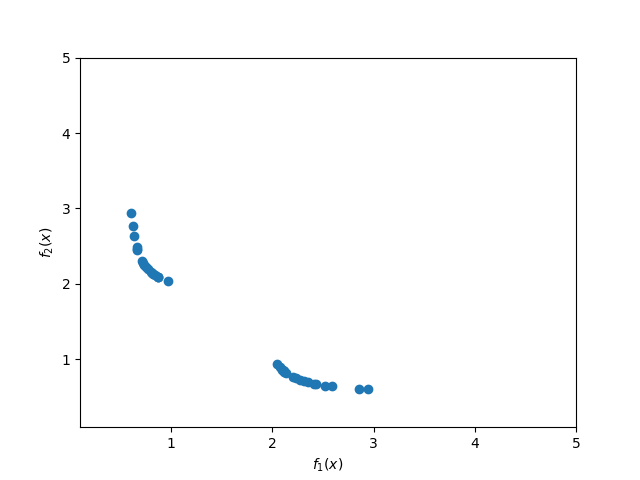}}
	\caption{Illustration of the clustering of the weighted averages in non-convex settings. For an idea of the shape of the true front see Figure~\ref{fig:with_penalization}.}
	\label{fig:A_interactionF}
\end{figure}

Figure~\ref{fig:A_interactionF} shows the approximation of the Pareto front and of the  efficient set that is obtained when using the algorithm with adaptive weights. We observe that the weighted means form clusters along the convex part of the Pareto front while the non-convex part is not recovered.

To overcome this issue, we propose a penalization strategy that avoids clusters along the Pareto front and tends to a uniform distribution of points along the Pareto front. Towards this end, we consider the potential leading to the interaction forces above, which is given by \eqref{eq:potential}. As we want to penalize clusters in the objective space, we use the distance in the objective space leading to the penalization term
\begin{equation}\label{eq:penalization}
	\fcluster (X_t^{k,j},v) = \sum_{\ell=1, \ell \ne k}^{K} R_c \, \e^{-|f(X_t^{k,j}) - f(\vtl)|/r_c}, \quad R_c, r_c > 0.
\end{equation}
By construction, the penalty term is smaller the farther away the objective of a particle is from the objective of the weighted means of the other swarms. The penalization term is added to the scalarization leading to a new cost given by
\[
f_{\mu^k,p}(X_t^{k,j}) = f_{\mu^k}(X_t^{k,j}) + \frac{\beta}{\alpha}\, \fcluster (X_t^{k,j},v), \qquad \alpha>0, \beta \ge 0.
\] 
Here $\beta$ allows us to balance the objective and the penalization.
Using the properties of the exponential function, we can rewrite the weighted average from \eqref{eq:weightedaverage} as
\begin{equation}\label{eq:weighted_average_penalization}
	{\vtk} = \frac{\sum_{j=1}^{N_k} X_t^{k,j} \e^{-\alpha f_{\mu^k}(X_t^{k,j})} \e^{-\beta\, \fcluster (X_t^{k,j},v)} }{\sum_{j=1}^{N_k}  \e^{-\alpha \,f_{\mu^k}(X_t^{k,j})} \e^{-\beta \fcluster (X_t^{k,j},v})}.
\end{equation}
Clearly, for $\beta = 0$ we are in the previous case without penalization of clustering in the objective space. With increasing values of $\beta >0$ the impact of this penalization in the optimization process is increased.
A pseudo code can be found in Algorithm~\ref{alg:everything}. 

\subsection{Numerical results for adaptive weights and  penalization}
\begin{figure}[ht!]
	\centering
	\subfloat[\label{}]{\includegraphics[scale=0.38]{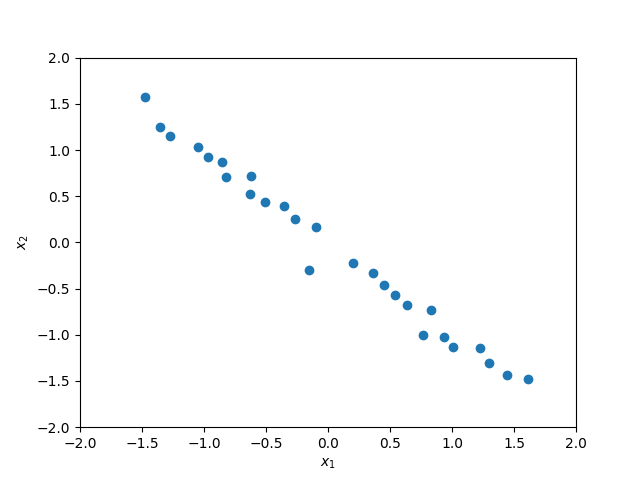}}
	\subfloat[\label{}]{\includegraphics[scale=0.38]{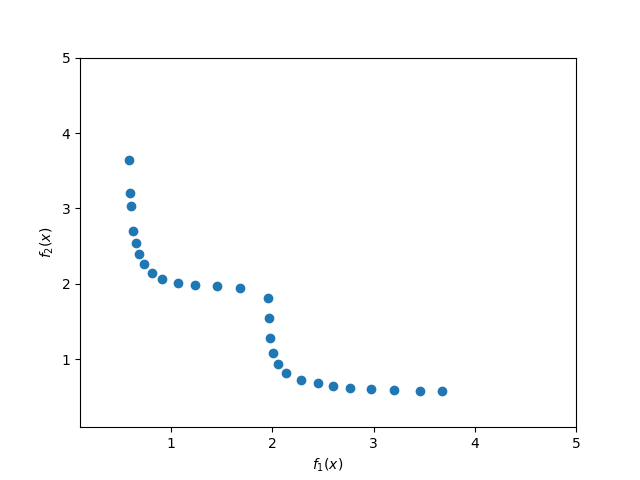}} 
	\caption{Illustration of the effect of penalization. On the left-hand side we see the approximation of the efficient set by the weighted means of the swarms obtained with Algorithm~\ref{alg:adaptive}. On the right-hand side the corresponding approximation of the Pareto front is shown.}
	\label{fig:with_penalization}
\end{figure}

The following results are obtained with the adaptive weight dynamics with additional penalization given in \eqref{eq:weighted_average_penalization}. We observe that the penalization prevents the clustering and the swarms are able to approximate also the non-convex part of the Pareto front. Moreover, the swarms spread along the convex parts leading to a better result in both the decision- and the  objective space, 
as is illustrated in Figure~\ref{fig:with_penalization}. The improvement compared to Figure~\ref{fig:A_interactionF} is obvious.

\begin{figure}
	\centering
	\subfloat[\label{}]{\includegraphics[scale=0.38]{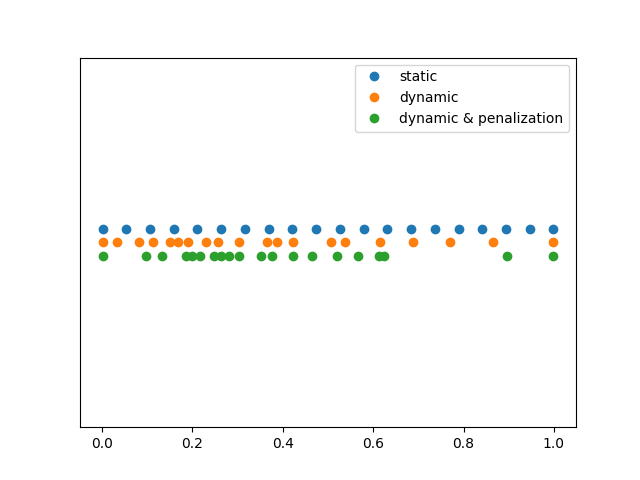}}
	\subfloat[\label{}]{\includegraphics[scale=0.38]{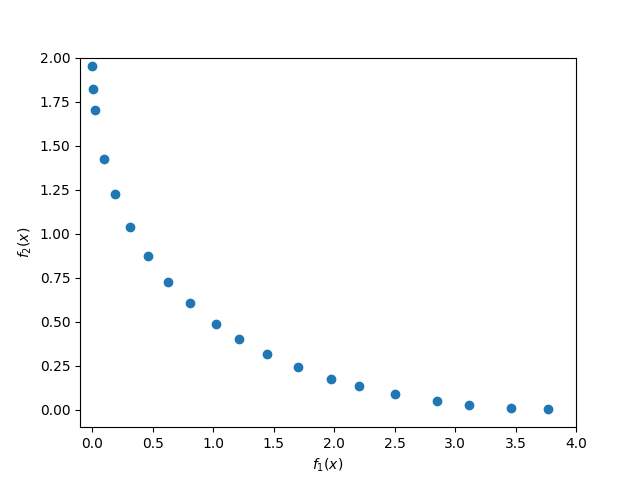}}
	\caption{Influence of penalization on the convex toy problem discussed in Section~\ref{sec:results_static_dynamic}.}
	\label{fig:convex_with_penalization}
\end{figure}

The influence of the penalization on the convex example of Section~\ref{sec:results_static_dynamic} is illustrated in Figure~\ref{fig:convex_with_penalization}. The approximation of the Pareto front obtained with penalization nicely covers the whole efficient set. In particular, the tails are better resolved compared to the approximation without penalization, cf.~Figure~\ref{fig:static_dynamics_weightsA}.

We emphasize that the penalization only affects the cost function of the problem. Therefore, all the analytical results obtained in Section~\ref{sec:simple} still apply. The only difference is that the weighted averages of the swarms approximate the minimizers of the penalized scalarizaed cost functions instead of the global minima of the (original) scalarized cost functions.

The theory of CBO provides us with a convergence result showing that the invariant state of the dynamics is a  consensus near the global minimizer of the cost function. We adapted the dynamics to obtain a scheme that provides a uniform approximation of the Pareto front all at once. So far, we only considered the weighed means of the swarms. Now, we want to use the additional information provided by the individuals. In fact, for our purpose it is better if the swarms do not collapse but cluster around the weighted mean, in order to obtain a local approximation of the Pareto front. 
We achieve this by implementing an idea from Consensus-based sampling (CBS) \cite{CarrilloHoffmannStuartVaes}. Technically speaking, we replace the anisotropic diffusion factor in Algorithm~\ref{alg:adaptive} by a factor that counteracts collapse of the swarm, see Algorithm~\ref{alg:everything}.

\begin{figure}[ht!]
	\centering
	\subfloat[\label{fig:diffusion_samplingA}]{\includegraphics[scale=0.38]{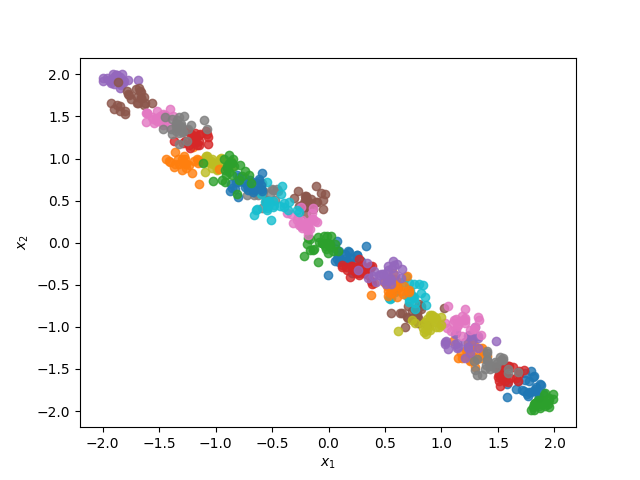}}
	\subfloat[\label{fig:diffusion_samplingB}]{\includegraphics[scale=0.38]{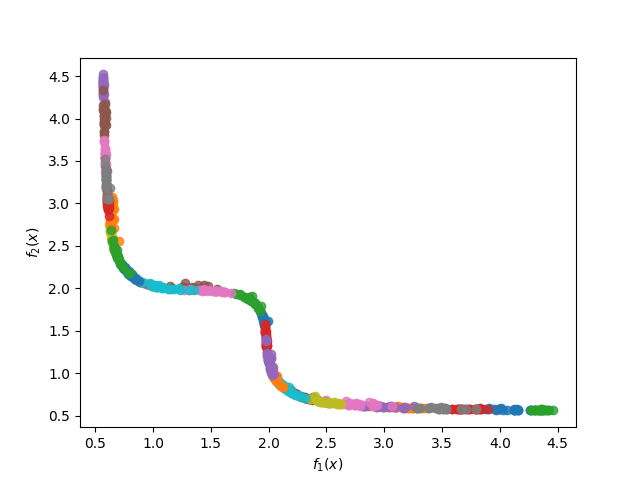}}
	\caption{Approximation of the Pareto front using the information of weighted means and all individuals. Individuals of the same swarm are plotted in the same color.}
	\label{fig:diffusion_sampling}
\end{figure}

Instead of the approximation consisting of the weighted means shown in Figure~\ref{fig:A_interactionF} and Figure~\ref{fig:with_penalization} we obtain now an approximation based on the weighted means and additional information of all individuals as shown in Figure~\ref{fig:diffusion_sampling}. The different colors of individuals illustrate the allocation to the respective swarms, plotted in the decision space  (Fig.~\ref{fig:diffusion_samplingA}) and in the objective space  (Fig.~\ref{fig:diffusion_samplingB}), respectively.

\begin{remark}
	Let us summarize the main ideas: in contrast to many evolutionary multiobjective optimization algorithms (EMO), like e.g., NSGA2 \cite{pymoo}, the proposed algorithms are not based on dominance tests of the swarms and/or particles. 
	Moreover, the adaptive weights automatically adjust the scalarization while exploring the objective functions. Finally, the diffusion term inspired by the sampling algorithm yields non-collapsing swarms that provide local information about the Pareto front near the diversely distributed weighted means.
\end{remark}

\section{Algorithms and parameter settings}
\label{sec:implementation}
In the following we provide the implementation details such as discretization, handling of box constraints and the parameter values for the simulations underlying the illustrations above.

\subsection{Discretization} \label{sec:discretization}
The stochastic differential equations modelling the multi-swarm dynamics are solved with the Euler-Maruyama scheme \cite{KloedenPlaten}. In each iteration particles are projected to the feasible set $\mathcal X$ (if necessary). Since we study only box-shaped feasible sets, the projection can be applied component-wise.
A pseudo code for the simple multi-swarm CBO  with fixed scalarization weights is given in Algorithm~\ref{alg:fixed}.

\begin{algorithm}
	\DontPrintSemicolon
	\caption{MSCBO 
		with fixed weights}\label{alg:fixed}
	\KwData{$K \in \mathbb N, N_k \in \mathbb N$, initial positions $\{X_0^{k,j}\}$ for $k=1,\dots,K$ and $j=1,\dots,N_k,$ scalarization weights $\lambda^k\in\Lambda_0$, time step $\tau>0$, diffusion coefficient $\sigma>0,$ terminal time $T> \tau,$ initial time $t=0$ and independent standard-normal distributed $W_t^k$}
	\KwResult{approximation of the Pareto front $(f_{\lambda^k}(X_T^{k,j}))_{k,j}$ and efficient set $(X_T^{k,j})_{k,j}$}
	
	\While{$t < T$}{
		\For {$k\gets1$ \KwTo $K$}{
			${\vtk} \gets \frac{\sum_{j=1}^{N_k} X_t^{k,j} \exp(-\alpha}f_{\lambda^k}(X_t^{k,j})){\sum_{j=1}^{N_k}  \exp(-\alpha f_{\lambda^k}(X_t^{k,j}))}$ \;
			$W_t^{k}\sim N(0,I)$\;
			$Y_{t+\tau}^{k} \gets X_t^{k} - \tau \, (X_t^{k}  - {\vtk}) + \sigma\, \sqrt{\tau}\, \diag(X_t^{k}  - {\vtk}) \, W_t^{k}$\;
			$X_{t+\tau}^{k} \gets \argmin\limits_{x \in \X} \|x - Y_t^{k}\|_2$ \;
		}
		$t \gets t+\tau$\;
	}
\end{algorithm}

The dynamic weight adaption proposed in Section~\ref{sec:adaptiveWeights} couples the dynamics of the different swarms. We proceed swarm by swarm and update the scalarization weight in each iteration explicitly. This leads to the pseudo code in Algorithm~\ref{alg:adaptive}.
\begin{algorithm}
	\DontPrintSemicolon
	\caption{MSCBO with adaptive weights}\label{alg:adaptive}
	\KwData{$K \in \mathbb N, N_k \in \mathbb N$, initial positions $\{X_0^{k,j}\}$ for $k=1,\dots,K$ and $j=1,\dots,N_k,$ scalarization weights $\mu^k = \ln(\lambda^k),$ time step $\tau>0$, diffusion coefficient $\sigma>0,$ terminal time $T> \tau,$ initial time $t=0$ and independent standard-normal distributed $W_t^k$}
	\KwResult{approximation of the Pareto front
		$(f_{\mu^k}(X_T^{k,j}))_{k,j}$ 
		and efficient set $(X_T^{k,j})_{k,j}$}
	\While{$t < T$}{
		\For {$k\gets1$ \KwTo $K$}{
			${\vtk} \gets \frac{\sum_{j=1}^{N_k} X_t^{k,j} \exp(-\alpha f_{\mu^k}(X_t^{k,j}))}{\sum_{j=1}^{N_k}  \exp(-\alpha f_{\mu^k}(X_t^{k,j}))}$ or 
			$\frac{\sum_{j=1}^{N_k} X_t^{k,j} \exp(-\alpha f_{\mu^k}(X_t^{k,j}) \, \exp(-\beta \fcluster (X_t^{k,j},v))}{\sum_{j=1}^{N_k}  \exp(-\alpha f_{\mu^k}(X_t^{k,j})) \, \exp(-\beta \fcluster (X_t^{k,j},v)) }$\;
			$W_t^{k}\sim N(0,I)$\;
			$Y_{t+\tau}^{k} \gets X_t^{k} - \tau (X_t^{k}  - {\vtk}) + \sigma \sqrt{\tau} \diag(X_t^{k}  - {\vtk}) W_t^{k}$\;
			$X_{t+\tau}^{k} \gets \argmin\limits_{x \in \X} \|x - Y_t^{k}\|_2$ \;
			$\mu_{t+\tau}^k \gets \mu_k^k - \frac{\tau}{K} \sum_{\ell=1,\ell \ne k}^K \mathcal K(X_t^k, X_t^\ell, \mu^k, \mu^\ell)$
		}
		$t \gets t+\tau$\;
	}
\end{algorithm}
The third modelling step implements the penalization strategy. Note that this affects only the cost function, hence leading to the adapted weighted mean (c.f.\ \eqref{eq:weighted_average_penalization} above)
\[
{\vtk} \gets \frac{\sum_{j=1}^{N_k} X_t^{k,j} \, \e^{-\alpha f_{\mu^k}(X_t^{k,j})}\, \e^{-\beta\,\fcluster (X_t^{k,j},v)}}{\sum_{j=1}^{N_k}  \e^{-\alpha \, f_{\mu^k}(X_t^{k,j})}\, \e^{-\beta \fcluster (X_t^{k,j},v)} }
\]

To avoid collapsing swarms we replace the diffusion term with the sampling noise leading to
\[
Y_{t+\tau}^{k} \gets X_t^{k} - \tau (X_t^{k}  - {\vtk}) + \sigma \sqrt{\tau |X_t^{k}  - {\vtk}|} \, W_t^{k}
\]

The pseudo code with all modelling features is summarized in Algorithm~\ref{alg:everything}.
\begin{algorithm}
	\DontPrintSemicolon
	\caption{MSCBO with adaptive weights, penalization and sampling noise}\label{alg:everything}
	\KwData{$K \in \mathbb N, N_k \in \mathbb N$, initial positions $\{X_0^{k,j}\}$ for $k=1,\dots,K$ and $j=1,\dots,N_k,$ scalarization weights $\mu^k = \ln(\lambda^k),$ time step $\tau>0$, diffusion coefficient $\sigma>0,$ terminal time $T> \tau,$ initial time $t=0$ and independent standard-normal distributed $W_t^k$}
	\KwResult{approximation of the Pareto front $(f_{\mu^k}(X_T^{k,j}))_{k,j}$ and efficient set $(X_T^{k,j})_{k,j}$}
	\While{$t < T$}{
		\For {$k\gets1$ \KwTo $K$}{
			${\vtk} \gets  
			\frac{\sum_{j=1}^{N_k} X_t^{k,j} \exp(-\alpha f_{\mu^k}(X_t^{k,j}))\, \exp(-\beta \fcluster (X_t^{k,j},v))}{\sum_{j=1}^{N_k}  \exp(-\alpha f_{\mu^k}(X_t^{k,j})) \,\exp(-\beta \fcluster (X_t^{k,j},v)) }$\;
			$W_t^{k}\sim N(0,I)$\;
			$Y_{t+\tau}^{k} \gets X_t^{k} - \tau (X_t^{k}  - {\vtk}) + \sigma \sqrt{\tau |X_t^{k}  - {\vtk}|}\, W_t^{k}$\;
			$X_{t+\tau}^{k} \gets \argmin\limits_{x \in \X} \|x - Y_t^{k}\|_2$ \;
			$\mu_{t+\tau}^k \gets \mu_k^k - \frac{\tau}{K} \sum_{\ell=1,\ell \ne k}^K \mathcal K(X_t^k, X_t^\ell, \mu^k, \mu^\ell)$
		}
		$t \gets t+\tau$\;
	}
\end{algorithm}

\subsection{Test problems}\label{sec:testproblems}
The illustrative example and also the numerical comparisons in the following sections use the following test problems.
\paragraph{Schaffer1 \cite{Schaffer}} \(f(x) = \big( (x-2)^2, \; 0.5\,x^2\big),\)
biobjective test problem, convex Pareto front, $\mathcal X = [0,2], \mathcal Y \subset [0,4]\times[0,2].$

\paragraph{Dent \cite{Witting}} $f(x) = (f_1(x_1,x_2), \; f_2(x_1,x_2))$ with
{\small\begin{align*} 
		f_1(x_1,x_2) &= \frac{1}{2} \left( \sqrt{1 + (x_1\!+\!x_2)^2} + \sqrt{1 + (x_1\! -\! x_2)^2} + x_1\! -\! x_2 \right) + 0.85 \e^{-(x_1 - x_2)^2},\\
		f_2(x_1,x_2) &= \frac{1}{2}\left( \sqrt{1 + (x_1\! +\! x_2)^2} + \sqrt{1 + (x_1\! -\! x_2)^2}  -x_1\! +\! x_2 \right) + 0.85 \e^{-(x_1 - x_2)^2}.
\end{align*}}
Biobjective test problem with a dent, $\mathcal X = [-2,2]^2, \mathcal Y \subset [0,5]^2.$

\paragraph{Schaffer2} $f(x) = \bigl(f_1(x),\; f_2(x)\bigr)$ with
{\small\begin{equation*}
		f_1(x) = \begin{cases} 
			-x, &\text{if } x\leq 1, \\ 
			x-2, &\text{if } 1\leq x\leq 3, \\
			4-x, &\text{if } 3\leq x\leq 4, \\
			x-4, &\text{if } 4\leq x \\
		\end{cases}, \qquad 
		f_2(x) = (x-5)^2.
\end{equation*}}
Biobjective test problem, discontinuous Pareto front, $\mathcal X = [-5,10], \mathcal Y \subset [0,1]\times[0,16].$

\paragraph{Three \cite{art:BOLT21}} 
$f(x_1,x_2) = \bigl(f_1(x_1,x_2),\; f_2(x_1,x_2),\; f_3(x_1,x_2)\bigr)$
with
{\small
	\begin{align*}
		f_1(x_1,x_2) &= 2(x_1 - 1)^2 + 2(x_1 - 1)(x_2 - 1) + 4(x_2- 1)^2, \\
		f_2(x_1,x_2) &= (x_1 - 2)^2 + 4(x_1 - 2)(x_2 - 3) + 8(x_2 - 3)^2, \\
		f_3(x_1,x_2) &= 4x_1^2 + 2x_1x_2 + x_2^2.
\end{align*}}
Three objectives, convex Pareto front, $\mathcal X = [-0.5,3.5], \mathcal Y \subset [0,25]\times[0,80]\times[0,50].$

\subsubsection*{Parameters for illustrative examples}\label{sec:parameters}
In the previous sections we illustrated the modelling ideas with the help of a convex optimization problem (\textsf{Schaffer1}) in Figures~\ref{fig:static_dynamics_weights} and~\ref{fig:convex_with_penalization}, and with a non-convex problem (dent) in Figures~\ref{fig:A_interactionF},~\ref{fig:with_penalization} and~\ref{fig:diffusion_sampling}. Initial data is sampled independently from the uniform distribution on $\mathcal X.$ All potentials are chosen to be purely repulsive, therefore $A$ and $A_f$ are set to zero, for completeness we set $a$ and $a_f$ to $1$. In a post-processing we eliminate dominated individuals, we therefore use a numerical offset $\epsilon_\text{dom} = 10^{-5}.$ For the test problems with $\mathcal Y \subset \mathbb R^2$ we initialize the scalarization weights equidistantly on $[0.001,0.999]^2$. 

For the three-objective test case, initial scalarization weights are randomly generated with Algorithm~2 in \cite{uniformSimplex}. For simplicity, we set $N_k = \bar N$ for all $k=1,\dots K$ for all simulations. For the illustrations in Figure~\ref{fig:static_dynamics_weights} and Figure~\ref{fig:diffusion_sampling} we use $K=20$ swarms with $\bar N=50$ individuals each. For all other biobjective test problems we set $K=30$ and $\bar N=20.$ For the test problem with three objectives we set $K=50$ and $\bar N=20.$ All other parameters are reported in Table~\ref{tab:parameters}.

\begin{table}[ht!]
	\centering
	\[
	\begin{array}{*{4}{r@{\extracolsep{0.5ex}}l@{\extracolsep{4ex}}}}
		\toprule
		\tau &= 0.1 & T &=  5 &\lambda &= 1 & \sigma &= 0.1 \\
		R &=  0.001 &  r &=  0.01 &  R_f &=  0.0001  &  r_f  &=  1  \\
		R_c &=  1  &  r_c  &=  0.1  &  \alpha  &= 100 &  \beta  &=  0  \text{ or }  10 \\\bottomrule
	\end{array}
	\]
	\caption{Algorithmic parameters used in the numerical tests.}
	\label{tab:parameters}
\end{table}

\section{Comparison to other population-based algorithms}\label{sec:comparison}
This section is devoted to a comparison to other multi-objective optimization methods. In particular, we address the well-known non-dominated Sorting Genetic Algorihm (NSGA2) \cite{NSGA2} and the recently introduced single swarm Consensus-based optimization approach for multi-objective problems \cite{BorghiHertyPareschi}. 

\medskip

\subsection{Computational effort}
The computational complexity of population-based multi-objective optimization algorithms lies on the one hand in the evaluation of the objective functions and on the other hand in computing the dominance relations between the generations or agents. The improvement from NSGA \cite{NSGA} to NSGA2 \cite{NSGA2} is mainly obtained by an efficient computation of the dominance relationship of the best $N$ individuals. Note that we have refrained from using the explicit multiobjective dominance information in the algorithms proposed and therefore save the complexity of $\mathcal O(pN^2)$ which is governed by the non-dominated sorting algorithms (see e.g.\ \cite{lang22space}).

The single swarm multi-objective algorithm (SSCBO) proposed in \cite{BorghiHertyPareschi} is stabilized by a greedy approach which ensures that individuals only move if the attempted move leads to a better position. Computationally this is cheaper then computing the dominance relationship in NSGA2, still it requires $N$ comparisons in each iteration. Our multi-swarm CBO approach MSCBO does not require these comparisons and therefore saves $\mathcal O(pN)$ as compared to SSCBO. 

\subsection{Numerical results}
To obtain comparable results, we chose the parameters in
the following simulations such that all algorithms use the same number of function evaluations. The time steps of SSCBO and MSCBO are chosen equally. We set the diffusion parameter of SSCBO to $10$ as reported in \cite{BorghiHertyPareschi}. Moreover, we activate the greedy strategy. The reference solution is given by NSGA2 with the same number of iterations. In the following this will be $T/\tau = 50.$

The comparison is based on different performance indicators: generational distance (GD), inverted generational distance (IGD) and hypervolume (HV) from the \texttt{pymoo} package are used for the comparison (see, e.g., \cite{riquelme15performance,zitzler03performance} for a detailed description of performance indicators). GD gives us an indication on the precision of the approximation of the Pareto front. IGD detects clustering or sparse regions along the front and HV gives us a flavour of how good the algorithms perform compared to NSGA2. As reference solution for GD and IGD we use the corresponding NSGA2 solution. The reference points for the hypervolume indicator are the upper bounds of the Cartesian intervals containing $\mathcal Y$ given in Section~\ref{sec:testproblems}. These choices allow us to analyse the performance without knowledge about the true Pareto front of the problems. As mentioned above we do not compute dominance relations in each iteration, but we evaluate the number of non-dominated individuals (NI) at the end of the simulations. 

\subsubsection{Biobjective tests}
\begin{figure}[ht!]
	\centering
	\subfloat[\label{fig:PFschaffer1A}]{\includegraphics[scale=0.38]{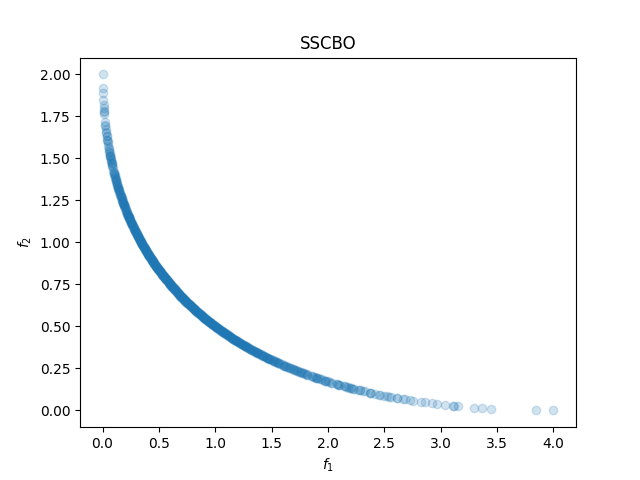}}
	\subfloat[\label{fig:PFschaffer1B}]{\includegraphics[scale=0.38]{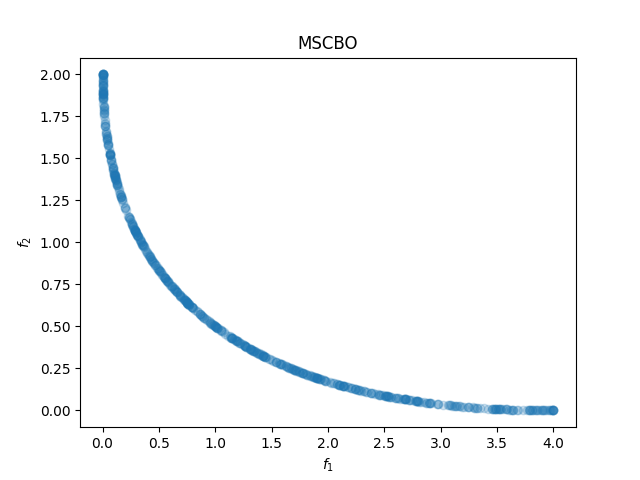}}
	\caption{Approximation of the Pareto fronts for \textsf{Schaffer1}.}
	\label{fig:PFschaffer1}
\end{figure}

We begin our tests with the convex test problem \textsf{Schaffer1} with the results reported in Table~\ref{tab:schaffer1}. In the outputs all individuals of all algorithms are non-dominated. The dominated hypervolumes of the three approximations are very similar. The GD values of MSCBO and SSCBO are very close as well, indicating that the approximations are precise. However, we see a difference in the IGD values. Further simulations and Figure~\ref{fig:PFschaffer1} suggest that the approximation points of SSCBO tend to accumulate in the region of the knee of the Pareto front, therefore the distances of IGD at the tails of the front are higher leading to a higher overall value. 

\begin{table}[ht!]
	\centering
	\begin{tabular}{l|c|c|c|c}
		$K=30, \bar N=20$        & GD & IGD & HV & NI\\ \midrule
		MSCBO & 0.00258 & 0.00396 & 6.65835 & 630\\
		SSCBO & 0.00227 & 0.01150 & 6.65686 & 630 \\
		NSGA2 & -- & -- & 6.66112 &  630 \\
	\end{tabular}
	\caption{Simulation results for \textsf{Schaffer1}.}
	\label{tab:schaffer1}
\end{table}

\begin{figure}[ht!]
	\centering
	\subfloat[\label{fig:PFdentA}]{\includegraphics[scale=0.38]{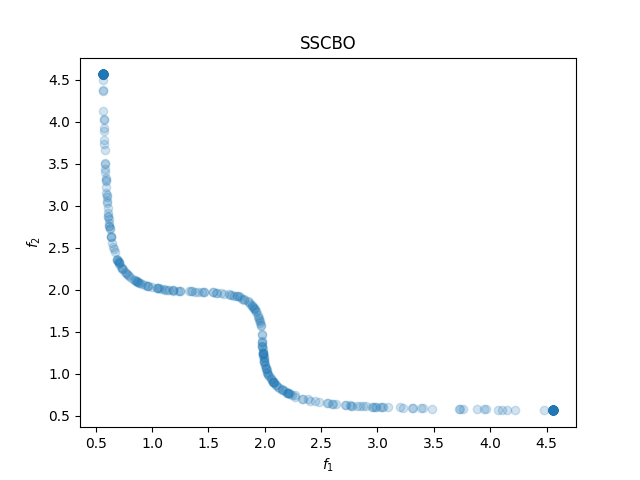}}
	\subfloat[\label{fig:PFdentB}]{\includegraphics[scale=0.38]{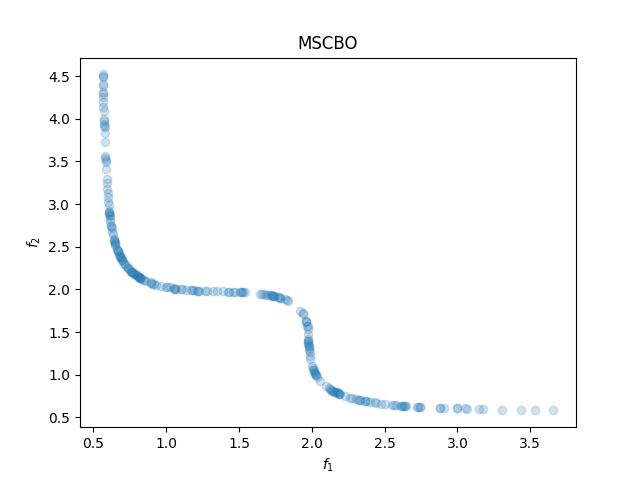}}
	\caption{Approximation of the Pareto fronts for problem \textsf{Dent}.}
	\label{fig:PFdent}
\end{figure}

The second test yields approximations of the Pareto front of the test problem \textsf{Dent}. It turns out that many of the individuals of MSCBO are dominated, also individuals of SSCBO are dominated. The GD values are very small for both methods. The IGD is better for SSCBO. Further investigations suggest that the approximation of the dent region of the Pareto front is sparser for MSCBO which leads to the higher values of IGD, see Figure~\ref{fig:PFdent}. Although MSCBO and SSCBO terminate with less non-dominated individuals, the dominated hypervolumes of all three methods are in good agreement.

\begin{table}[ht!]
	\centering
	\begin{tabular}{l|c|c|c|c}
		$K=30, \bar N=20$        & GD & IGD & HV & NI\\ \midrule
		MSCBO & 0.00388 & 0.06476 & 9.49007 & 230 \\
		SSCBO & 0.00244 & 0.01761 & 9.49547 & 375 \\
		NSGA2 & -- & -- & 9.51303 & 630 \\
	\end{tabular}
	\caption{Simulation results for problem \textsf{Dent}.}
	\label{tab:dent}
\end{table}

\begin{figure}[ht!]
	\centering
	\subfloat[\label{fig:PFschaffer2A}]{\includegraphics[scale=0.38]{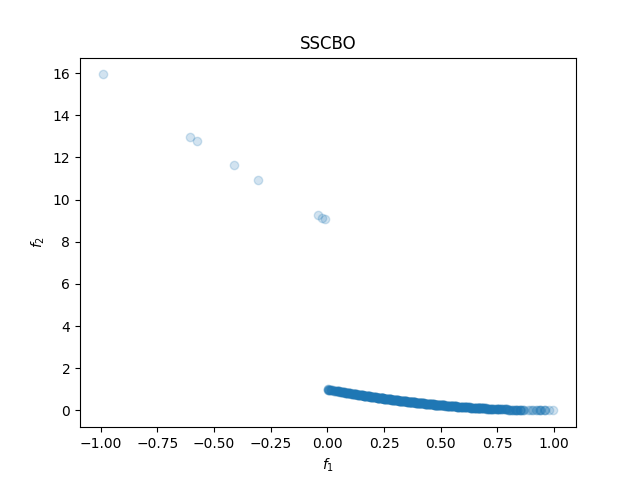}}
	\subfloat[\label{fig:PFschaffer2B}]{\includegraphics[scale=0.38]{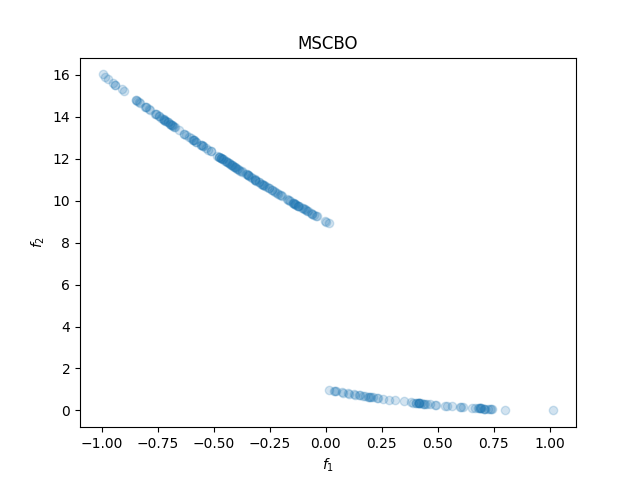}}
	\caption{Approximation of the Pareto fronts for \textsf{Schaffer2}.}
	\label{fig:PFschaffer2}
\end{figure}

The third test is concerned with the \textsf{Schaffer2} problem which has a discontinuous Pareto front. The precision of the approximation of the front indicated by the GD value is again very good. The IGD values differ and simulations showed that the upper part of the Pareto front is not very well-approximated by SSCBO leading to this IGD behavior. This is visualized in  Figure~\ref{fig:PFschaffer2} as well. Also in the values of the hypervolume indicator we see a difference between the SSCBO approximation. Although many of the MSCBO individuals are dominated, the dominated hypervolume and the IGD values are promising.

\begin{table}[ht!]
	\centering
	\begin{tabular}{l|c|c|c|c}
		$K=30, \bar N=20$    & GD & IGD & HV & NI\\ \midrule
		MSCBO &  0.00518 & 0.02456& 19.15976 &221 \\
		SSCBO & 0.00214 & 0.28919 & 18.37781 & 615\\
		NSGA2 & -- & -- & 19.31961 & 630 \\
	\end{tabular}
	\caption{Simulation results for \textsf{Schaffer2}.}
	\label{tab:schaffer2}
\end{table}

\subsubsection{Three-objective test}
The results of the benchmark with three objectives are shown in Table~\ref{tab:three}. It jumps to the eye that many individuals of SSCBO are dominated. This is not expected as the front is convex. For all tested indicators MSCBO outperforms SSCBO. There seems to be a performance drop from two to three objectives in SSCBO. The approximation determined by MSCBO is reasonable. Note that we do expect the magnitude of the indicators to be higher due to the additional objective function.

\begin{table}[ht!]
	\centering
	\begin{tabular}{l|c|c|c|c}
		$K=50, \bar N=20$   & GD & IGD & HV & NI\\ \midrule
		MSCBO & 0.54036 & 2.5566 & 88359 & 1041\\
		SSCBO & 0.60967 & 4.8240 & 87408 & 669\\
		NSGA2 & -- & -- & 89288 & 1050\\
	\end{tabular}
	\caption{Simulation results for three-objective problem. }
	\label{tab:three}
\end{table}

\section{Conclusion and outlook}
\label{sec:conclusion}
We propose a versatile population-based algorithm for the approximation of the Pareto front and of the efficient set for multi-objective optimization problems that is based on the Consensus-based optimization or sampling method, respectively. In the case of fixed scalarization weights, many analytical results of CBO can be easily generalized to this setting. For adaptive weights we proved that the methods yields a diverse approximation of the front. Moreover, with prevention of collapsing swarms, we retain local information of the Pareto front close to the swarm means. The method is computationally cheap as no dominance relations need to be computed. The numerical results are promising and motivate for further investigations and applications of the method. In fact, the algorithm is competitive with widely-used NSGA2 and the recently proposed single-swarm CBO for multi-objective problems. 

\section*{Acknowledgement}
We thank the authors of \cite{BorghiHertyPareschi} for providing their code of the single swarm multi-objective CBO. This facilitated the comparison of the methods significantly. We acknowledge the open source code of NSGA2 and the performance indicators in the project \texttt{pymoo} \cite{pymoo}.

\bibliographystyle{alpha}
\bibliography{biblio}

\end{document}